\documentclass[smallextended]{svjour3}       

\usepackage{alg}
\usepackage{amsfonts}
\usepackage{amsmath}
\usepackage{amssymb}
\usepackage{bm}
\usepackage{color}
\usepackage{graphicx}
\usepackage{subfigure}

\graphicspath{{../Figures/}}

\def\lf{\left}
\def\rt{\right}
\newcommand{\norm}[1]{\lVert#1\rVert}

\newcommand\RR{\mathbb{R}}
\newcommand\dt{\Delta t}
\newcommand\dx{\Delta x}
\newcommand\dy{\Delta y}
\newcommand\ds{\Delta s}
\newcommand\s{\bm{s}}
\newcommand\x{\bm{x}}
\newcommand\X{\bm{X}}
\newcommand\F{\bm{F}}
\newcommand\f{\bm{f}}
\newcommand\e{\bm{e}}
\renewcommand\u{\bm{u}}
\newcommand\Af{A_{\rm f}}   
\newcommand\w{\bm{w}}
\newcommand\Amg{\mathcal{A}}

\newcommand{\change}[1]{#1}

\begin{document}

\title{Geometric multigrid for an implicit-time immersed boundary method
\thanks{This work was supported in part by University of
    California Office of The President (UCOP) grant
    09-LR-03-116724-GUYR to RG and BP, National Science Foundation
    (NSF) grant DMS 1160438 to RG, and by American Heart Association
    grant 10SDG4320049 and NSF grants DMS 1016554 and OCI 1047734 to
    BG.}
}

\author{Robert D.\ Guy \and Bobby Philip \and Boyce E.\ Griffith}

\institute{Robert D.\ Guy \at
              Department of Mathematics, University of California, Davis, CA \\
              Tel.: +1-530-754-9201\\
              Fax: +1-530-752-6635 \\
              \email{guy@math.ucdavis.edu}           
           \and
           Bobby Philip \at 
           Oak Ridge National Laboratory, Oak Ridge, TN \\
           \email{philipb@ornl.gov}       
           \and 
           Boyce E.\ Griffith \at
              Leon H.\ Charney Division of Cardiology, Department of
              Medicine, New York University School of Medicine, New York, NY \\
              \email{boyce.griffith@nyumc.org}
}

\date{Received: date / Accepted: date}

\maketitle

\begin{abstract}
  The immersed boundary (IB) method is an approach to fluid-structure
  interaction that uses Lagrangian variables to describe the
  deformations and resulting forces of the structure and Eulerian
  variables to describe the motion and forces of the fluid.  Explicit
  time stepping schemes for the IB method require solvers only for
  Eulerian equations, for which fast Cartesian grid solution methods
  are available.  Such methods are relatively straightforward to
  develop and are widely used in practice but often require very small
  time steps to maintain stability.  Implicit-time IB methods permit
  the stable use of large time steps, but efficient implementations of
  such methods require significantly more complex solvers that
  effectively treat both Lagrangian and Eulerian variables
  simultaneously.  Several different approaches to solving the coupled
  Lagrangian-Eulerian equations have been proposed, but a complete
  understanding of this problem is still emerging.  This paper
  presents a geometric multigrid method for an implicit-time
  discretization of the IB equations.  This multigrid scheme uses a
  generalization of box relaxation that is shown to handle problems in
  which the physical stiffness of the structure is very large.
  Numerical examples are provided to illustrate the effectiveness and
  efficiency of the algorithms described herein.  These tests show
  that using multigrid as a preconditioner for a Krylov method yields
  improvements in both robustness and efficiency as compared to using
  multigrid as a solver.  They also demonstrate that with a time step
  100--1000 times larger than that permitted by an explicit IB method,
  the multigrid-preconditioned implicit IB method is approximately
  50--200 times more efficient than the explicit method.



\end{abstract}

%

%
%
\section{Introduction}
\label{intro:sec}

The immersed boundary (IB) method \cite{CambridgeJournals:165651} was
introduced by Peskin \cite{Peskin1977220} to solve problems of
fluid-structure interaction in which an elastic structure is immersed
in a viscous incompressible fluid.  This method was developed to
simulate the dynamics of heart valves, but it has subsequently been
applied to diverse problems in biofluid dynamics, and it is finding
increasing use in other engineering problems
\cite{MITTAL:2005:IBREVIEW}.  The IB formulation of such problems uses
an Eulerian description of the momentum, viscosity, and
incompressibility of the fluid-structure system, and it uses a
Lagrangian description of the deformation of the immersed structure
and forces generated by these deformations.  The Eulerian equations
are approximated on a Cartesian grid, and the Lagrangian equations are
approximated on a curvilinear mesh.  Interaction between Eulerian and
Lagrangian variables is through integral equations with delta function
kernels.  When discretized, the IB method uses a regularized version
of the delta function to mediate Lagrangian-Eulerian coupling.  A key
feature of the method is that it does not require conforming
discretizations of the fluid and structure; instead, the curvilinear
mesh is free to cut through the background Cartesian grid in an
arbitrary manner.  Consequently, IB simulations do not require dynamic
grid generation, even for problems involving very large structural
deformations.

Typical implementations of the IB method adopt a fractional step
approach to time stepping.  In the simplest version of such a scheme,
the Eulerian velocity and pressure fields are updated for a fixed
configuration of the immersed structure, and then the position of the
Lagrangian structure is updated from the newly computed velocity
field.  This approach effectively decouples the Eulerian and
Lagrangian equations, and solvers are needed only for the Eulerian
equations (i.e., the incompressible Stokes or Navier-Stokes
equations), for which fast Cartesian grid solution methods are
available.  However, because this fractional step approach yields an
explicit time stepping method for the structural dynamics problem,
maintaining stability requires time steps that are small enough to
resolve all of the elastic modes of the discrete equations.  In many
applications, these elastic time scales are well below the physical
time scales of interest.  Even for relatively simple elasticity
models, the largest stable time step size scales like $\dt =
O(\ds^2)$, in which $\ds$ is the Lagrangian mesh spacing.  For
problems involving bending-resistant elastic elements, the largest
stable time step scales like $\dt = O(\ds^4)$.  Consequently,
high-resolution IB simulations can require extremely large numbers of
time steps, making it challenging to perform simulations over long
time scales.

Much effort has been devoted both to understanding and to alleviating
the severe time step restriction of fractional step IB methods
\cite{Stockie199941,Newren2007702,gong:2008:ibanalysis}.  One
approach is to develop implicit or semi-implicit time stepping schemes
that allow for time steps that do not resolve all of the elastic modes
of the discrete system; however, despite decades of work, such schemes
are still not widely used in practice.  The solution methods used in
early implicit IB methods were not efficient and were not competitive
with explicit methods \cite{tu:1361}, and some semi-implicit methods
intended to allow for large time steps still suffered from significant
time step restrictions \cite{MAYO:1993:IMPLICIT,lee:832}.  Newren et
al.\ \cite{Newren2007702} analyzed the origin of instability in
semi-implicit IB methods using energy arguments, and they gave
sufficient conditions to achieve unconditionally stability in the
sense that the total energy is bounded independent of time step size.
An important result by Newren et al.\ \cite{Newren2007702} is that it
is not necessary to employ a \emph{fully} implicit time discretization
to achieve unconditional stability, but the stable time stepping
schemes proposed therein do simultaneously solve for both the Eulerian
velocity field and the Lagrangian structural configuration.  However, as
indicated by the early experience with implicit IB methods
developing efficient solvers for the coupled equations is challenging.

More recently, a number of stable semi-implicit \cite{Newren20082290,Hou20088968,Hou20089138,Ceniceros20097137} and fully implicit
\cite{Mori20082049,Le20098427} IB methods have been developed.  The
efficiency of these methods is generally competitive with explicit
methods, and in some special cases, these implicit schemes can be
faster than explicit methods by several orders of magnitude.  Many
implicit methods use a Schur complement approach to reduce the coupled
Lagrangian-Eulerian equations to purely Lagrangian equations
\cite{Mori20082049,Ceniceros20097137,ceniceros:3d_implicit_ib:2011}.
These methods achieve a substantial speed-up over explicit methods
when there are relatively few Lagrangian mesh nodes
\cite{Ceniceros20097137}.  In addition, some methods require that the
boundaries be smooth, closed curves \cite{Hou20088968,Hou20089138}.
An open question is whether there exist robust, general-purpose
implicit methods that are more efficient than explicit methods, or
whether specialized methods must be developed for specific problems.

Newren et al.\ \cite{Newren20082290} explored the use of
unpreconditioned Krylov methods for solving the linearized IB
equations in the context of different test problems.  They found that
the relative efficiency of the implicit methods depended on the
problem, and unpreconditioned Krylov methods were generally at least
comparable in speed to explicit methods.  These results suggest that
with appropriate preconditioning, this approach will offer a
significant improvement over explicit methods.  One way to achieve
generally applicable and robust implicit methods is through the
development of robust preconditioners for the linearized equations.
This is the approach we take in this paper.

In previous work \cite{GUY:2012:IMPLICITIB}, we developed a multigrid
method for a model problem related to implicit time discretizations of
the IB equations.  This model problem ignored the inertial terms and
the incompressibility constraint.  The multigrid solver introduced in
this earlier work was more efficient than explicit-time methods for
the model problem, but the increase in efficiency was not large for
the very stiff problems for which implicit time stepping methods are
most clearly needed.  When used as a preconditioner for a Krylov
solver, however, the method was very efficient, even for very stiff
problems.

In this paper, we extend the methods developed for the model IB
equations \cite{GUY:2012:IMPLICITIB} to problems of incompressible
flow.  Specifically, we consider a version of the IB method for the
steady incompressible Stokes equations.  (The extension of the method
to the unsteady Stokes equations, or to the full Navier-Stokes
equations, is straightforward but is not considered here.)  Unlike
most other work on developing efficient solvers for implicit IB
methods, here we focus on a formulation of the problem in which we
effectively eliminate the \emph{structural} degrees of freedom by a
Schur complement approach.  The system that we solve is therefore
defined only on the background Cartesian grid.  As in earlier work
\cite{GUY:2012:IMPLICITIB}, we take advantage of the structure of the
Cartesian grid to develop geometric multigrid methods for the linear
systems of our implicit-time discretization.  The key contributions of
this paper are the development of generalized box-relaxation (also
known as Vanka) smoothers for this formulation of the IB equations,
and the extension of these box-relaxation smoothers to larger
collections of grid cells, as needed to obtain good performance for
problems in which the elastic structure is extremely stiff.  We
perform numerical tests that demonstrate the performance of these
algorithms, and we show that with these solvers, the implicit scheme
has the potential to be significantly more efficient than a similar
explicit IB method.

%
%
\section{Immersed Boundary Equations}
\label{ibeqn:sec}

\subsection{Continuum equations}

Let $\x\in\Omega$ denote fixed physical coordinates, with
$\Omega\subset\RR^2$ denoting the physical domain.  Let $\s\in\Gamma$
denote material coordinates attached to the immersed structure, with
$\Gamma\subset\RR^2$ denoting the Lagrangian coordinate
domain.\footnote{We remark that the name \emph{immersed boundary
    method} suggests that the elastic structure is a thin interface
  (i.e., an object of codimension one with respect to the fluid).
  While this is the case in many applications of the IB method, this
  formulation applies equally well to immersed structures that have
  nonzero thickness.  We restrict our tests to two spatial dimensions
  and to structures of nonzero thickness.  The extension to three
  spatial dimensions is straightforward, and in the concluding
  discussion, we comment on the differences between thick and thin
  structures.}  The physical location of material point $\s$ at time
$t$ is given by $\X(\s,t)\in\Omega$.  (In general, we use lowercase
letters for quantities expressed in Eulerian coordinates and uppercase
letters for quantities expressed in Lagrangian coordinates.)  In the
absence of other loading, the forces generated by the deformations of
the structure drive the motion of the fluid.  We assume that the
immersed structure is neutrally buoyant, so that all of the boundary
force is transmitted to the fluid.  The equations we consider in this
paper are
\begin{gather}
  \Delta \u(\x,t) - \nabla p(\x,t) + \bm{f}(\x,t) = 0, \\
  \nabla\cdot \u(\x,t) = 0, \\
  \bm{f}(\x,t) = \int_{\Gamma} \bm{F}(\s,t)\,\delta\!\lf(\x-\X(\s,t)\rt)\,{\mathrm d}\s, \label{spread:eq} \\
  \frac{\partial \X(\s,t)}{\partial t} = \bm{U}(\s,t)
         = \int_{\Omega} \u(\x,t) \, \delta\!\lf(\x-\X(\s,t)\rt)\,{\mathrm d}\x,  \label{interp:eq}
\end{gather}
in which $\u(\x,t) = (u(\x,t),v(\x,t))$ is the velocity field of the
fluid-structure system, $p(\x,t)$ is the pressure, $\bm{f}(\x,t)$ is
the Eulerian elastic force density generated by the immersed
structure, and $\bm{F}(\s,t)$ is the Lagrangian elastic force density
generated by the immersed structure.  The first two equations are the
incompressible Stokes equations, which here describe the motion of a
fluid-structure system in which the influence of inertia is
negligible.  The last two equations describe the coupling between the
Eulerian and Lagrangian frames.  The integral operator in
\eqref{spread:eq} that determines the Eulerian force density from the
Lagrangian force density is called the \emph{spreading} operator,
which we denote by $S[\X]$.  The \emph{interpolation} operator that
transfers the velocity to the structure is the adjoint of the
spreading operator.  Using this notation, equations \eqref{spread:eq}
and \eqref{interp:eq} can be compactly expressed as $\bm{f} =
S[\X]\,\bm{F}$ and $\partial \X / \partial t = \bm{U} =
S[\X]^{*}\,\u$, respectively.

\begin{figure}[t]
  \centering
  \includegraphics[width=0.45\textwidth]{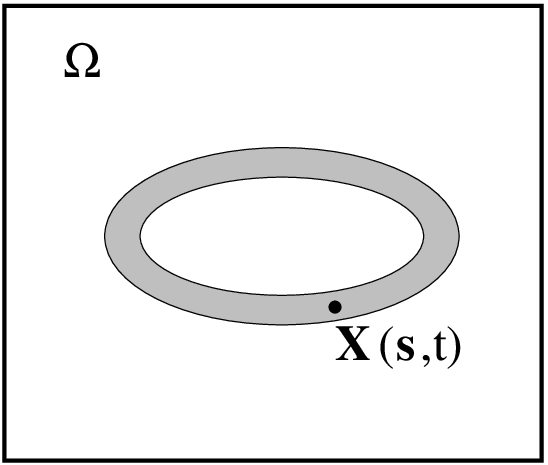}
  \caption{The physical domain $\Omega$ contains the immersed elastic
    structure.  The position of each material point $\s$ at time $t$
    is given by $\X(\s,t)\in\Omega$.}
  \label{domain2:fig}
\end{figure}

A constitutive law for the immersed elastic material is needed to
complete the description of the system.  Herein, we consider
structures that consist of a collection of linear elastic fibers under
tension.  We choose the Lagrangian coordinate system so that for
$\s=(s_{1},s_{2})$, $s_{1}$ is a parametric coordinate along each
fiber, and $s_{2}$ is constant on each fiber.  Let $\bm{\tau}$ be the
unit vector tangent to the fiber direction, which is given by
$\bm{\tau} = \partial\X/\partial s_{1} / \norm{\partial\X/\partial
  s_{1}}$.  The tension in each fiber is taken to be $T = \gamma
\norm{\partial\X/\partial s_{1}}$, in which $\gamma$ is a constant
that characterizes the elastic stiffness of the fiber.  Under these
assumptions, the Lagrangian force density is
\begin{equation}
  \F(\s,t) = \frac{\partial}{\partial s_{1}}\lf(T\bm{\tau}\rt)
           = \gamma \frac{\partial^{2}\X}{\partial s_{1}^{2}}.
           \label{ibforce:eq}
\end{equation}
This constitutive law corresponds to an elastic shell that is composed
of a continuum of circumferential elastic fibers
\cite{Boffi20082210,griffith:adaptib:2007,Griffith200575}.  It is also
equivalent to a version of an incompressible neo-Hookean elastic
material.  \change{Note that in this formulation, the overall
  structural response is actually viscoelastic, because the structural
  stress tensor includes a viscous term that is identical to that of
  the fluid.  Immersed boundary methods have been developed which
  include a separate structure viscosity
  \cite{FAI:2013:VarVisc,Huang20092650,STRYCHALSKI:2012:VEIB}.  We
  take the structure viscosity to be equal to the fluid viscosity here
  for simplicity; however, as discussed in Section
  \ref{discussion:sec}, this is not a fundamental limitation of our
  method.}

\subsection{Spatial discretization}

\begin{figure}[hb]
  \centering
  \includegraphics[width=0.45\textwidth]{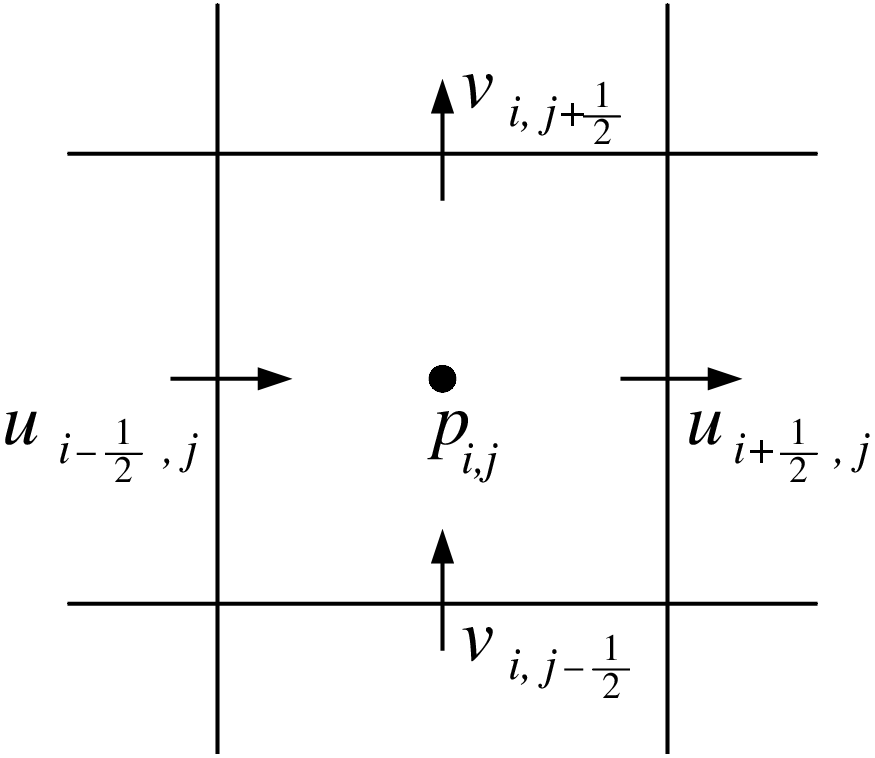}
  \caption{Staggered grid discretization in which the velocity field
    $\u=(u,v)$ and Eulerian body force $\f=(f_1,f_2)$ are approximated
    at the centers of the cell edges, and in which the pressure $p$ is
    approximated at the cell centers.}
  \label{macproj:fig}
\end{figure}

The physical domain $\Omega$ is taken to be rectangular and in our
computations is discretized by a uniform Cartesian grid with square
cells of width $\dx = \dy = h$.  We use a staggered-grid
discretization of the incompressible Stokes equations in which the
components of the velocity and Eulerian body force are approximated at
the centers of the cell edges to which that component is normal, and
in which the pressure is approximated at the cell centers; see Figure
\ref{macproj:fig}.  The Laplacian, gradient, and divergence operators
are discretized using standard second-order finite differences, and
the corresponding discrete operators are denoted by $\Delta_{h}$,
$\nabla_{h}\cdot\mbox{}$, and $\nabla_{h}$, respectively.

The immersed structure is discretized using a fiber-aligned mesh with
nodes that are equally spaced in the Lagrangian coordinate system with
spacing $\ds_1 = \ds_2 = \ds$ is each direction.  The physical
location of Lagrangian node $\s_{k,l}$ is denoted by $\X_{k,l}$.  The
second derivative operator that appears in the boundary force
\eqref{ibforce:eq} is discretized using standard three-point centered
differencing.  The Lagrangian elastic force density at node $\s_{k,l}$
is
\begin{equation}
   \F_{k,l}=
   \gamma \frac{\X_{k-1,l} - 2\X_{k,l} + \X_{k+1,l}}{\ds^2} = \gamma \lf(\Af \X\rt)_{k,l},
\end{equation}
in which we use $\Af$ to denote the discrete force-generation
operator.

To approximate the Lagrangian-Eulerian interaction equations, we use a
two-dimensional regularized delta function $\delta_h(\x)$ that is the
tensor product of two one-dimensional regularized delta functions, \change{$\delta^{1}_{h}(r)$}, so that for $\x = (x,y)$ and $\X = (X,Y)$,
\begin{equation}
  \delta_h(\x-\X) = \delta^{1}_{h}(x-X) \, \delta^{1}_{h}(y-Y).
\end{equation}
In this work, we use
\begin{equation}
  \delta^{1}_{h}(r) =
  \begin{cases}
    \frac{1}{4h}\lf(1 + \cos\lf(\frac{\pi r}{2}\rt) \rt) & \text{if } r < 2h, \\
      0                   & \text{otherwise.}
  \end{cases}
\end{equation}
The discretized spreading operator $S_h[\X]$ is defined for
$\f = (f_1,f_2)$ and for $\F = (F_1,F_2)$ by
\begin{gather}
  \lf(f_1\rt)_{i-1/2,j}=
  \sum_{k,l}\lf(F_1\rt)_{k,l} \, \delta_{h}(\x_{i-1/2,j}-\X_{k,l}) \, \ds^{2}, \\
  \lf(f_2\rt)_{i,j-1/2}=
  \sum_{k,l}\lf
(F_2\rt)_{k,l} \, \delta_{h}(\x_{i,j-1/2}-\X_{k,l}) \, \ds^{2},
\end{gather}
in which $\x_{i-1/2,j}$ and $\x_{i,j-1/2}$ denote the positions of the
centers of the edges of the grid cells, where the velocity and force
components are approximated.  Similarly, the discrete interpolation
operator, which is the adjoint of the discrete spreading operator, is
defined by
\begin{gather}
  \lf(\frac{\partial X}{\partial t}\rt)_{k,l} = U_{k,l} =
  \sum_{i,j}u_{i-1/2,j} \, \delta_{h}(\x_{i-1/2,j}-\X_{k,l}) \, h^{2}, \\
  \lf(\frac{\partial Y}{\partial t}\rt)_{k,l} = V_{k,l} =
  \sum_{i,j}v_{i,j-1/2} \, \delta_{h}(\x_{i,j-1/2}-\X_{k,l}) \, h^{2}.
\end{gather}


\subsection{Temporal discretizations}

\subsubsection{Explicit-time method}

Typical implementations of the IB method use a fractional time
stepping approach to solve the equations of motion.  In the simplest
version of such a scheme, the fluid velocity and pressure are updated
while keeping the position of the structure fixed, and then the
structural position is updated using the newly computed velocity.  We
refer to this method as the explicit-time method.  For the model
equations considered herein, the explicit-time method advances the
solution variables from time $t^{n}=n\dt$ to time $t^{n+1}=(n+1)\dt$
via
\begin{gather}
  \Delta_{h} \u^{n+1} - \nabla_{h} p^{n+1} + S_h^{n} \F^{n} = 0, \label{explicit1:eq}\\
  \nabla_{h}\cdot\u^{n+1} = 0, \\
  \F^{n} = \gamma \Af X^{n} \\
  \X^{n+1} = \X^{n} + \dt \lf(S_h^{n}\rt)^{*} \u^{n+1}, \label{explicit2:eq}
\end{gather}
in which $S_h^{n} = S_h\lf[\X^n\rt]$.  Notice that the explicit-time
method requires the solution of only the incompressible Stokes system.

\subsubsection{Implicit-time method}

The implicit-time method is similar to the ex\-plic\-it-time method,
except that we now use backward Euler to update the structural
position, and we now compute the structural forces using the newly
computed positions.  The implicit method advances the solution
variables via
\begin{gather}
  \Delta_{h} \u^{n+1} - \nabla_{h} p^{n+1} + S_h^{n} \F^{n+1} = 0, \label{model_mom:eq} \\
  \nabla_{h}\cdot\u^{n+1} = 0, \\
  \F^{n+1} = \gamma \Af \X^{n+1}, \label{model_F:eq} \\
  \X^{n+1} = \X^{n} + \dt \lf(S_h^{n}\rt)^{*} \u^{n+1}. \label{model_X:eq}
\end{gather}
Notice that in this time stepping scheme, the structural positions
used to define the spreading and interpolation operators are lagged in
time.  As shown by Newren et al.\ \cite{Newren2007702}, this scheme is
unconditionally stable, despite the fact that the positions of the
spreading and interpolation operators are treated explicitly rather
than implicitly.  This is quite fortuitous; if it were necessary to
treat the spreading and interpolation operators implicitly, then we
would be faced with a \emph{nonlinear} system of equations.

We use \eqref{model_F:eq} and \eqref{model_X:eq} to eliminate the
unknown force $\F^{n+1}$ from \eqref{model_mom:eq} to yield a system
in which the only unknowns are the velocity $\u^{n+1}$ and pressure
$p^{n+1}$,
\begin{gather}
   \lf(\Delta_{h}+\alpha S_h^{n}\Af \lf(S_h^{n}\rt)^{*} \rt)\u^{n+1} -\nabla_{h} p + S_h^{n} \Af \X^{n} = 0  \label{model1:eq}\\
  \nabla_{h}\cdot \u^{n+1} = 0,  \label{model2:eq}
\end{gather}
with $\alpha = \dt\gamma$.

To advance the full system in time, we first solve equations
\eqref{model1:eq}\---\eqref{model2:eq} for the velocity and pressure,
and we then use equation \eqref{model_X:eq} to update the position of
the structure.  The advantages of reducing the full system
\eqref{model_mom:eq}\---\eqref{model_X:eq} to
\eqref{model1:eq}\---\eqref{model2:eq} are that the only unknowns in
\eqref{model1:eq}\---\eqref{model2:eq} are the Eulerian velocity and
pressure, and that the equations are defined on a structured grid.
The elimination of the Lagrangian unknowns facilitates the development
of geometric multigrid methods for the IB equations.

%
%
\section{Multigrid}
\label{mg:sec}

We provide a brief sketch of geometric multigrid methods, focusing on
details specific to our application.  For a detailed description, the
reader is referred to Refs.~\cite{TROTTENBERG:BOOK,Briggs00}.

Let $\Omega_h$ represent the discretized physical domain with
Cartesian grid spacing $h$.  The linear system
\eqref{model1:eq}--\eqref{model2:eq} on $\Omega_h$ can be written as
\begin{equation}
   \lf[ \begin{array}{cc}
      \lf(\Delta_{h}+\alpha S_h\Af S_h^{*} \rt)  &  -\nabla_{h} \\
      \nabla_{h}\cdot\mbox{}                     &    0
     \end{array} \rt]
  \lf[ \begin{array}{l} \u^{n+1} \\ p^{n+1} \end{array}\rt]
  =
  \lf[\begin{array}{c} -S_h\Af \X^{n} \\ 0\end{array}\rt],
\label{matrixeq:eq}
\end{equation}
which we denote by
\begin{equation}
  \Amg_h\bm{w}_h = \bm{b}_h. \label{eqn-linear-system}
\end{equation}
To simplify the notation, we set $S_h \equiv S_h^{n}$ and $S_h^{*}
\equiv (S_h^{n})^{*}$.  Notice, however, that $S_h$ and $S_h^{*}$ are
generally time-dependent discrete operators, as is $\Amg_h$.

Let $I_{2h \leftarrow h}$ denote the operator that \emph{restricts}
solution data to $\Omega_{2h}$ from $\Omega_h$, let $I_{h \leftarrow
  2h}$ denote the operator that \emph{prolongs} solution data to
$\Omega_h$ from $\Omega_{2h}$, and let $\mathcal{A}_{2h}$ denote the
coarse-grid operator defined on $\Omega_{2h}$.  (Definitions for each
of these operators are provided below.)  The smoothers used in this
work, which are a key aspect of the overall solution algorithm, are
specified below in the context of specific numerical examples.  We
specifically consider a generalization of a standard box-relaxation
smoother in Section \ref{box_relaxation:sec}, and an extension of this
approach to ``big'' boxes in Sections \ref{big-box_smoothing:sec} and
\ref{timedeptest:sec}.  A geometric multigrid V-cycle for
\eqref{eqn-linear-system} is given by Algorithm \ref{alg:mgcycle}.
\begin{algorithm}[h]
\caption{$\w_h \longleftarrow MGV(\w_h, \bm{b}_h, \Omega_h, \nu_1, \nu_2)$\label{alg:mgcycle}}
\begin{algtab}
  \algif {$\Omega_h$ is the coarsest level}
  Solve the coarse-grid equation: $\w_h \longleftarrow (\Amg_h)^{-1}\bm{b}_h$ \\
  \algelse
  Perform $\nu_1$ \emph{presmoothing} sweeps for $\Amg_h\w_h=\bm{b}_h$ on $\Omega_h$ using initial guess $\w_h$\\
  Compute the residual on $\Omega_h$ and restrict it from $\Omega_h$ to $\Omega_{2h}$: $\bm{r}_{2h}\,\:\longleftarrow I_{2h \leftarrow h}\lf(\bm{b}_h-\Amg_h\w_h\rt)$ \\
  Compute an approximate solution to the error equation $\Amg_h\bm{e}_h=\bm{r}_h$ on $\Omega_{2h}$: $\bm{e}_{2h}\longleftarrow MGV(\bm{0}, \bm{r}_{2h}, \Omega_{2h}, \nu_1, \nu_2)$\\
  Prolong the coarse-grid correction from $\Omega_{2h}$ to $\Omega_h$ and update the solution on $\Omega_h$: $\w_h\,\:\longleftarrow \w_h+I_{h \leftarrow 2h} \e_{2h}$\\
  Perform $\nu_2$ \emph{postsmoothing} sweeps for $\Amg_h\w_h=\bm{b}_h$ on $\Omega_h$ using initial guess $\w_h$\\
\end{algtab}
\end{algorithm}

Multigrid methods are intended to work on all components of the error
in a given approximation to the solution of \eqref{eqn-linear-system}
by a combination of fine-grid relaxation (steps 4 and 8 in Algorithm
\ref{alg:mgcycle}) and coarse-grid correction (steps 5--7 in Algorithm
\ref{alg:mgcycle}).  In well-designed multigrid methods, fine-grid
relaxation and coarse-grid correction are complementary processes: the
errors that are not damped by the fine-grid relaxation are damped by
the coarse-grid correction, and vice versa.  If these processes are
not complementary (i.e., do not damp all error modes), then the method
will yield poor convergence rates; if these processes overlap (i.e.,
damp the same components of the error), then the method will provide
sub-optimal efficiency.


\subsection{Grid transfer operators: Restriction and prolongation}

We use standard geometric coarsening of the Cartesian grid in which a
hierarchy of successively coarser grids $\Omega_h$, $\Omega_{2h}$,
$\Omega_{4h}$, \ldots is generated.  Because approximations to the
components of the velocity and the pressure are all defined at
different spatial locations, different operators are required to
transfer these values between levels of the hierarchy of
discretizations.

For the pressure, we obtain coarse cell-centered values by averaging
the four overlying fine cell-centered values.  The stencil and
coefficients of the pressure restriction operator are given by
\begin{equation}
  R_{p} = \frac{1}{4}\lf[
    \begin{array}{ccc}
       1 &   & 1 \\
         & * &   \\
       1 &   & 1
    \end{array}
    \rt] ,
\end{equation}
in which ``$*$'' denotes the position of the coarse value.  To prolong
pressure data from coarser grids to finer grids, we use constant
prolongation, so that for each coarse grid cell, each overlying fine
grid cell takes the underlying coarse grid value.

Restriction of the $x$-components of the velocity ($u$) is done by
two-point averaging in the $y$-direction and full-weighting in the
$x$-direction.  The stencil and coefficients for the operator are
\begin{equation}
  R_{u} = \frac{1}{8}\lf[
    \begin{array}{ccc}
      1 & 2 & 1 \\
        & * &   \\
      1 & 2 & 1
\end{array}\rt].
\end{equation}
A similar procedure, but with two-point averaging in the $x$-direction
and full-weighting in the $y$-direction, is used to restrict the $y$
components of the velocity ($v$).  In each case, standard bilinear
interpolation is used to prolong components of the velocity from
coarser grids to finer grids.

We remark that these transfer operators are the standard ones for
staggered-grid discretizations of incompressible flow problems
\cite{TROTTENBERG:BOOK}, but other transfer operators could be used.  See
Niestegge et al.\ \cite{Niestegge1990} for a study of the performance
of several different combinations of interpolation and restriction
operators for the Stokes equations.  We experimented with different
combinations of operators, and we found these standard operators gave
the best overall efficiency.


\subsection{Coarse-grid operator}

Coarse-grid correction requires the formation of a coarse-grid
operator for each of the coarser levels of the grid hierarchy.  In
geometric multigrid, the two most common approaches are direct
re-discretization of the PDE on each grid level, and algebraic
construction via a Galerkin procedure.  In previous work on a model of
the IB method, we found that Galerkin coarsening was necessary
for convergence \cite{GUY:2012:IMPLICITIB}.  However, Galerkin
coarsening of the Stokes equations is expensive because the Galerkin
coarse-grid operators have large stencils.  Because re-discretization
works well for the Stokes equations alone (i.e., without the IB
elasticity operator $S_h \Af S_h^{*}$), we employ a hybrid approach:
We re-discretize the Stokes equations, and we use Galerkin coarsening
for the Eulerian elasticity operator $S_h \Af S_h^*$.  Specifically,
the coarse-grid operator is
\begin{equation}
 \Amg_{2h} =   \lf[ \begin{array}{cc}
      \lf(\Delta_{2h}+\alpha  I_{2h \leftarrow h} S_h\Af S_h^{*} I_{h \leftarrow 2h} \rt)  &  -\nabla_{2h} \\
      \nabla_{2h}\cdot\mbox{}                      &    0
     \end{array} \rt].
\end{equation}
Coarser versions of $S_h \Af S_h^{*}$, e.g., on $\Omega_{4h}$,
$\Omega_{8h}$, \ldots, are constructed recursively.


\subsection{Multigrid preconditioning}

As we have remarked, multigrid is a highly effective solver when
smoothing and coarse-grid correction work in a complementary manner to
eliminate all error modes, and with the smoothers used in this work,
the present algorithm achieves high efficiency for the Stokes problem.
As the stiffness of the immersed boundary increases in the
implicit-time method, however, the discrete operator becomes
increasingly less ``Stokes-like'' in the vicinity of the immersed
structure, and the performance of the multigrid algorithm suffers.

In our previous work on applying multigrid to a model of the IB method
\cite{GUY:2012:IMPLICITIB}, we found that multigrid alone was a poor
solver for large stiffnesses but that it performed very effectively as
a preconditioner for Krylov methods.  We follow the same approach
here, and in our numerical experiments, we explore the performance of
multigrid as both a solver and as a (right) preconditioner for GMRES
\cite{saad:856}.  For more details on the general use of multigrid
preconditioned Krylov methods, see Refs.\ \cite{TATEBE:MGPCG:1993,Oosterlee:1998:EPM:277965.277972}, and see Refs.\
\cite{Sussman199981,Wright:2008:ERM:1405171.1405187,Elman1996} for
specific examples from fluid mechanics.

%
%
\section{Test Problem Description}

We explore the performance of the multigrid method as a solver and as
a preconditioner for a range of elastic stiffnesses of the immersed
structure.  Except where otherwise noted, the physical domain $\Omega$
is the unit square $[0,1]^{2}$, and Dirichlet conditions are imposed
on the velocity along $\partial\Omega$ to yield lid-driven-cavity
flow.  Specifically, all components of the velocity are set to zero on
the boundary except on the top wall, where the tangential velocity is
$u(x,1)=(1-\cos(2\pi x))/2$.

In all cases, the immersed structure is the annulus with initial
positions
\begin{equation}
 \X(s_{1},s_{2}) =   \Bigl(x_\text{c} + (r+s_{2})\cos(s_{1}/r),
                          y_\text{c} +(r+s_{2})\sin(s_{1}/r) \Bigr),
\label{thick_coords:eq}
\end{equation}
in which $\x_\text{c} = (x_\text{c},y_\text{c}) = (0.5,0.5)$ is the
center of the annulus, which also generally corresponds to the center
of $\Omega$ in our tests, and in which $r=1/4$ is the inner radius of
the annulus.  The Lagrangian coordinate domain is
$(s_{1},s_{2})\in[0,2\pi r)\times[0,w]$, with $w=1/16$ indicating the
thickness of the annulus.  This domain is discretized using a regular
grid with $M_{1}$ points in the $s_{1}$ direction and $M_{2}$ points
in the $s_{2}$ direction.  We choose $M_{1}=19 N/8$ and $M_{2}=3N/32 +
1$, in which $N$ is the number of grid cells used to discretize one
direction in the Eulerian domain, so that $h = 1/N$.  We restrict $N$
to be a power of two, so that $M_{1}$ and $M_{2}$ are integers.  The
physical distance between adjacent Lagrangian nodes is approximately
$2/3$ of the Eulerian grid spacing $h$.

\change{A similar thick ring structure has been used in other immersed
  boundary benchmarking tests on accuracy \cite{Griffith200575},
  volume conservation \cite{GRIFFITH:2012:IBVOLCONSERVE}, and implicit
  time stepping \cite{Mori20082049}.  A notable difference between our
  test and past tests is that we drive a background flow though the
  boundary condition so that the physical time scale is set by the
  background flow, not by the stiffness of the structure.  For large
  elastic stiffness, the elastic time scale is well below the physical
  flow time scale, which is a characteristic of problems that benefit
  most from efficient implicit-time methods.}

\subsection{Characterizing the elastic stiffness}

The explicit-time method given by equations
\eqref{explicit1:eq}\---\eqref{explicit2:eq} is equivalent to the
forward Euler scheme applied to
\begin{equation}
  \frac{\partial \X}{\partial t} = \gamma S_h^{*}\mathcal{L}_h^{-1} S_h \Af \X,
  \label{explicitX:eq}
\end{equation}
in which $\mathcal{L}_h^{-1}$ is the operator that maps fluid forces
to the fluid velocity by solving the Stokes system, and  
 $\gamma$ is the stiffness of the elastic structure (see
\eqref{ibforce:eq}). The stability of
this scheme is determined by the single parameter
\begin{equation}
  \alpha = \gamma\dt.
\end{equation}
\change{The forward Euler method applied to
  \eqref{explicitX:eq} is stable provided 
\begin{equation}
  \lvert 1 + \alpha\lambda\rvert \leq 1,
  \label{stab_fe:eq}
\end{equation}
for all $\lambda$ which are eigenvalues of
$S_h^{*}\mathcal{L}_h^{-1}S_h \Af$.  Because
$S_h^{*}\mathcal{L}_h^{-1}S_h$ is symmetric and positive semidefinite
and $\Af$ is symmetric and negative semidefinite, the eigenvalues of
their product are all real and nonpositive.  Therefore the stability
condition \eqref{stab_fe:eq} is equivalent to to the condition 
\begin{equation}
  \alpha \rho \leq 2,
\end{equation}
where $\rho$ is the spectral radius of the matrix
$S_h^{*}\mathcal{L}_h^{-1} S_h \Af$.}
Let $\alpha_\text{exp}$ denote the maximum value
of $\alpha$ for which the explicit-time method is stable, which is
defined by
\begin{equation}
  \alpha_\text{exp} = \frac{2}{\rho},
\end{equation}
In Table
\ref{alpha_max_balloon:tab}, we report values of $\alpha_\text{exp}$
for different grid spacings.

\begin{table}[ht]
\centering
\caption{$\alpha_\text{exp}$ is the maximum value of the stiffness
$\alpha=\gamma\dt$ for which the explicit-time scheme is stable for
grid spacing $h$.}
\begin{tabular}{|r|r|r|r|r|}
\hline
$h$                 & $2^{-5}$ & $2^{-6}$ & $2^{-7}$ & $2^{-8}$ \\ \hline
$\alpha_\text{exp}$ & 6.09     & 3.93     & 2.82     & 2.28     \\ \hline
\end{tabular}
\label{alpha_max_balloon:tab}
\end{table}

%
%
\section{Box Relaxation}
\label{box_relaxation:sec}

Several different smoothers for the Stokes equations have been
developed.  Two large classes of smoothers are distributive smoothers
\cite{BrandtDinar79} and collective smoothers \cite{Vanka86}.
Distributive relaxation techniques, originally pioneered by Brandt and
Dinar \cite{BrandtDinar79}, involve a transformation of the equations
so that the individual velocity components and pressure are smoothed
separately \cite{BrandtDinar79,Linden89,Wittum90}.  Collective or box
relaxation, originally proposed by Vanka \cite{Vanka86}, involves
smoothing the velocity and the pressure simultaneously.  Oosterlee and
Washio \cite{Oosterlee08} provide a comparison of distributed and
collective smoothers for incompressible flow problems.

Distributed smoothers for the Stokes equations are straightforward to
implement because they involve smoothing only scalar problems;
however, extending distributive smoothers to the implicit IB equations
is challenging.  In particular, it is not clear whether it is possible
to transform the saddle point problem \eqref{matrixeq:eq} into a form
that permits the decoupled smoothing of the velocity and pressure.  In
this work, we instead employ box relaxation.  Box relaxation is
essentially a generalization of point smoothers like Jacobi or
Gauss-Seidel to multi-component systems, including saddle-point
systems.  The basic idea of box relaxation is to sweep over the grid
cells and, in each cell, to solve locally the discrete equations
restricted to that cell.  In the present context, a 5-by-5 system of
equations must be solved for each cell that involves the four velocity
components and the one pressure.  We order the boxes
lexicographically, and we update the unknowns box-by-box in a block
Gauss-Seidel-like manner.

\subsection{Solver performance}

\begin{figure}[ht]
  \centering
  \includegraphics[width=0.5\textwidth]{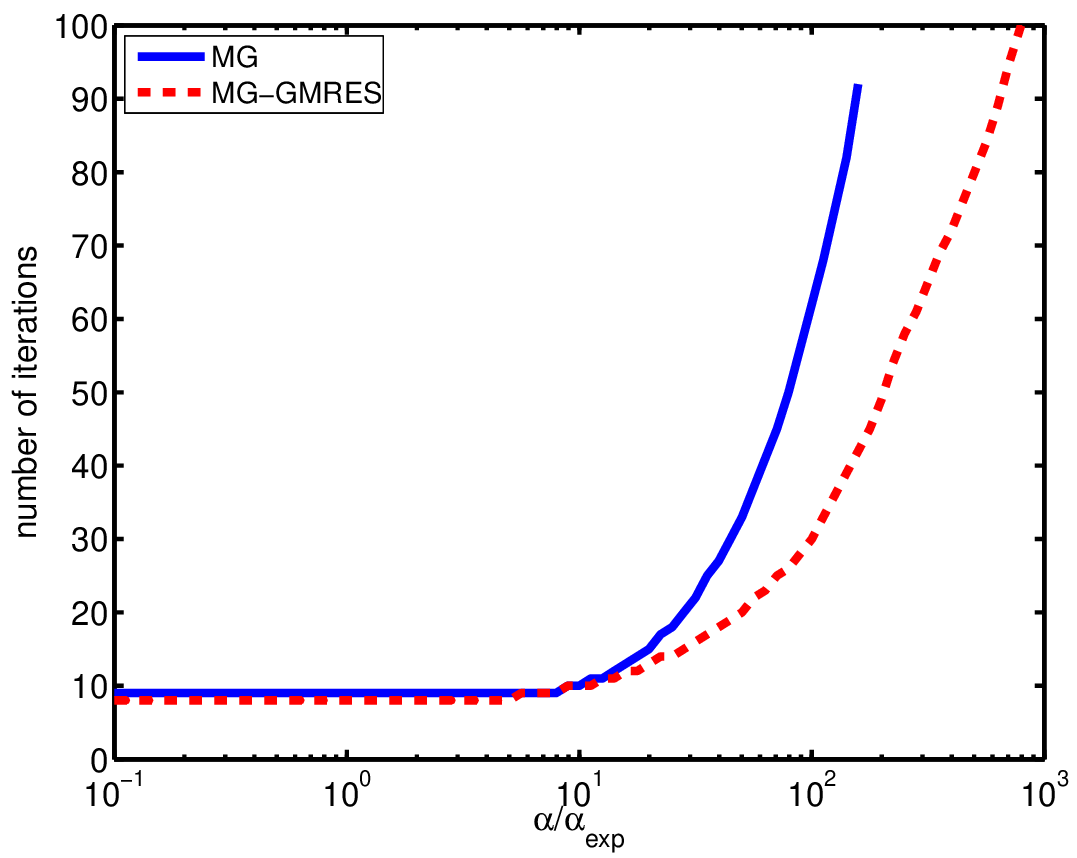}
  \caption{Iteration counts for the multigrid solver and for
    multigrid-preconditioned GMRES to reduce the residual by a factor
    of $10^{-6}$ as a function of the stiffness. The stiffness is
    scaled by the maximum stiffness from the explicit-time scheme,
    $\alpha_\text{exp}$.}
  \label{singlebox:fig}
\end{figure}

As an initial test of the performance of the multigrid method as a
solver and as a preconditioner, we consider an Eulerian grid with
$h=2^{-5}$ and the corresponding Lagrangian mesh, we set the initial
guess for the velocity and pressure to zero, and we compute the number
of iterations needed to reduce the residual by a factor of $10^{-6}$.
We use V-cycles with one presmoothing sweep and one postsmoothing
sweep ($\nu_1 = \nu_2 = 1$).  Figure \ref{singlebox:fig} shows the
resulting iteration counts as a function of $\alpha/\alpha_\text{exp}$
for both multigrid as a solver and as a preconditioner for GMRES.  The
ratio $\alpha/\alpha_\text{exp}$ may be interpreted as follows.  For a
given elastic stiffness, $\alpha/\alpha_\text{exp}$ represents the
size of the time step relative to the maximum allowed by the
explicit-time method.  We call this ratio the \emph{relative
  stiffness}.

The iteration count of the multigrid method is essentially constant (9
iterations) up to a relative stiffness of about $10$, at which point
the iteration count begins to increase rapidly.  For
$\alpha/\alpha_\text{exp}=100$, the multigrid solver takes 62
iterations to converge, and the multigrid solver fails to converge for
a relative stiffness of $\alpha/\alpha_\text{exp}\approx160$.

The iteration count for the MG-GMRES method is also essentially
constant (8 iterations) up to a relative stiffness of $10$.  After
this point, the iteration count begins to increase, but not as rapidly
as when we use multigrid as a solver.  For
$\alpha/\alpha_\text{exp}=100$, the iteration count is $30$, or about
half that of multigrid alone.  Unlike the stand-alone multigrid
solver, MG-GMRES does not appear ever to diverge, but we stopped after
100 iterations.  This maximum number of iterations was reached at a
relative stiffness of $\alpha/\alpha_\text{exp}\approx 800$.

If we assume that a V-cycle for the Stokes-IB system takes about the
same amount of work as a V-cycle applied to the Stokes equations
(i.e., without the IB elasticity operator), then we can use the
iteration counts to estimate the efficiency of this
approach.\footnote{In practice, the extent to which this assumption
  holds depends on the relative densities of the Eulerian and
  Lagrangian discretizations.}  For $\alpha/\alpha_\text{exp}=100$,
the multigrid algorithm takes 62 iterations.  To reach the same point
in time, an explicit method would require at least 100 time steps and
require 9 iterations per time step.  Therefore we estimate the
implicit method would be $100\cdot 9/62 \approx 14.5$ times more
efficient.  A similar estimate for MG-GMRES suggests that the implicit
method would be about 27 times more efficient at this value of
stiffness.  These are likely overestimates of the efficiency gain, and
we will return to more careful efficiency comparisons for a
time-dependent problem in Section \ref{timedeptest:sec}.

For very stiff problems MG-GMRES is more efficient than the
stand-alone multigrid solver.  For nonstiff problems (e.g.,
$\alpha/\alpha_\text{exp}<10$), stand-alone multigrid is more
efficient because one iteration of MG-GMRES is more expensive than one
iteration of multigrid.  However, MG-GMRES is the more robust solver;
it does not fail to converge as the stiffness increases.


\subsection{Spectrum of the multigrid operator}

To explore the relatively poor performance and ultimate failure of the
stand-alone multigrid solver at large stiffness, we explicitly
construct the multigrid iteration matrix and compute its eigenvalues.
We construct the multigrid iteration matrix using the procedure
outlined in Ref.\ \cite{TROTTENBERG:BOOK}.  The idea is as follows:
Let the values of the velocity and pressure be organized into a single
vector $\bm{w}=[u; v; p]$.  To generate the $k^\text{th}$ column of
the multigrid matrix, 
we set $\bm{w}_{j}=\delta_{jk}$, in
which $\delta_{jk}$ is the Kronecker delta, and we perform one
multigrid cycle.

\begin{figure}[h]
 \centering
   \subfigure[$\alpha/\alpha_\text{exp}=0,~ \rho\approx 0.24 $]
           {\includegraphics[width=0.42\textwidth]{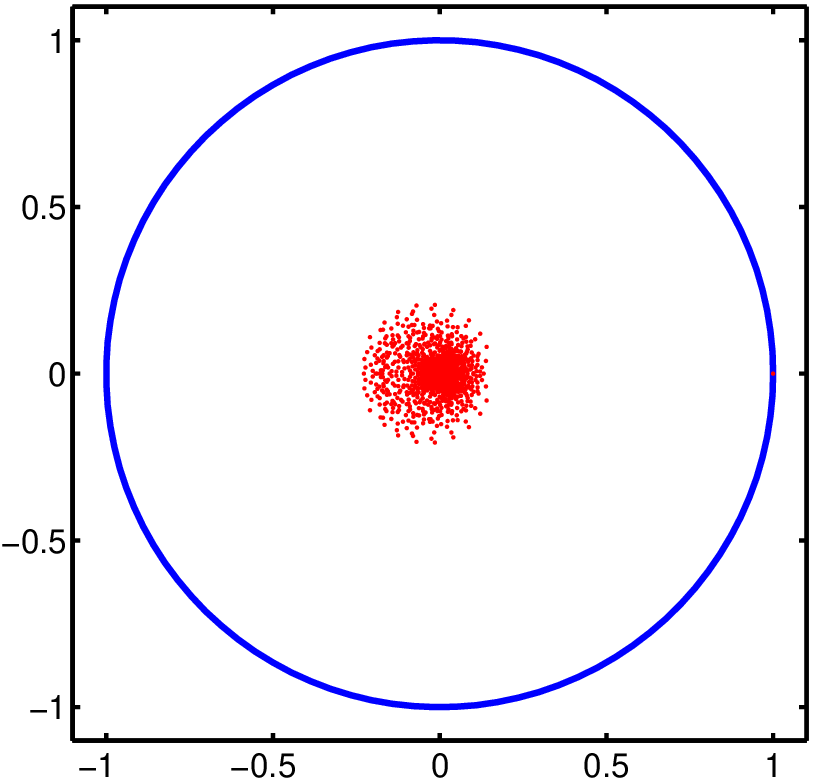}} \hfill
   \subfigure[$\alpha/\alpha_\text{exp}=1,~ \rho\approx 0.22$]
           {\includegraphics[width=0.42\textwidth]{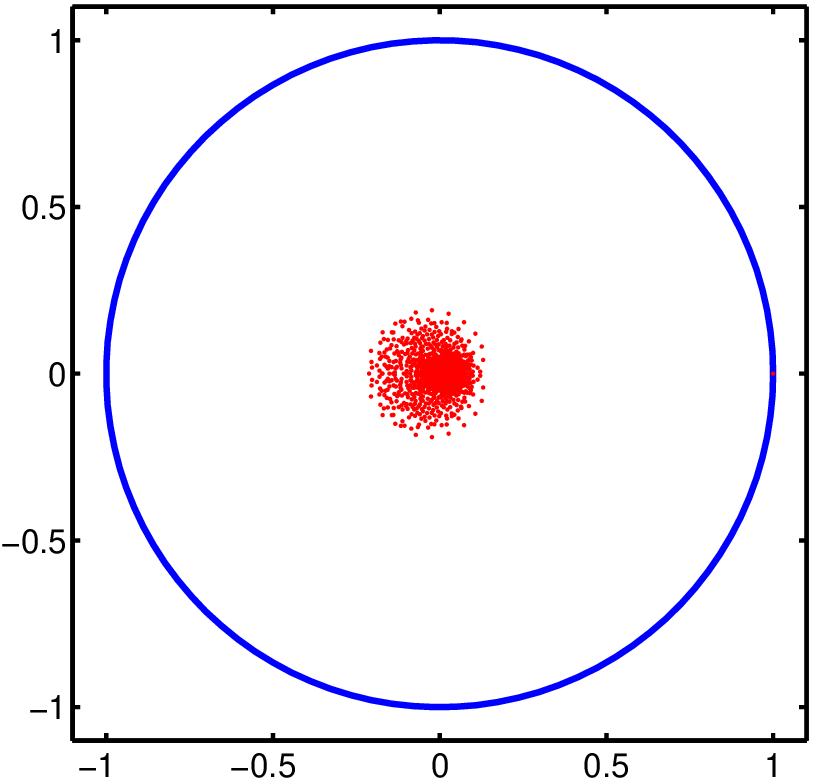}}
   \subfigure[$\alpha/\alpha_\text{exp}=10,~ \rho\approx  0.30$]
           {\includegraphics[width=0.42\textwidth]{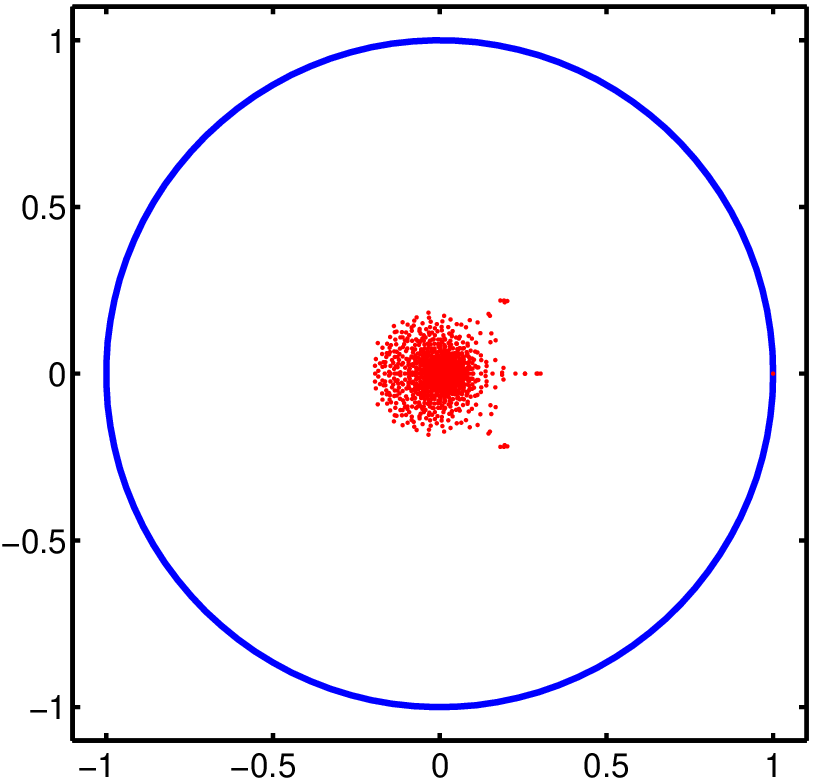}} \hfill
   \subfigure[$\alpha/\alpha_\text{exp}=100,~ \rho\approx  0.85$]
           {\includegraphics[width=0.42\textwidth]{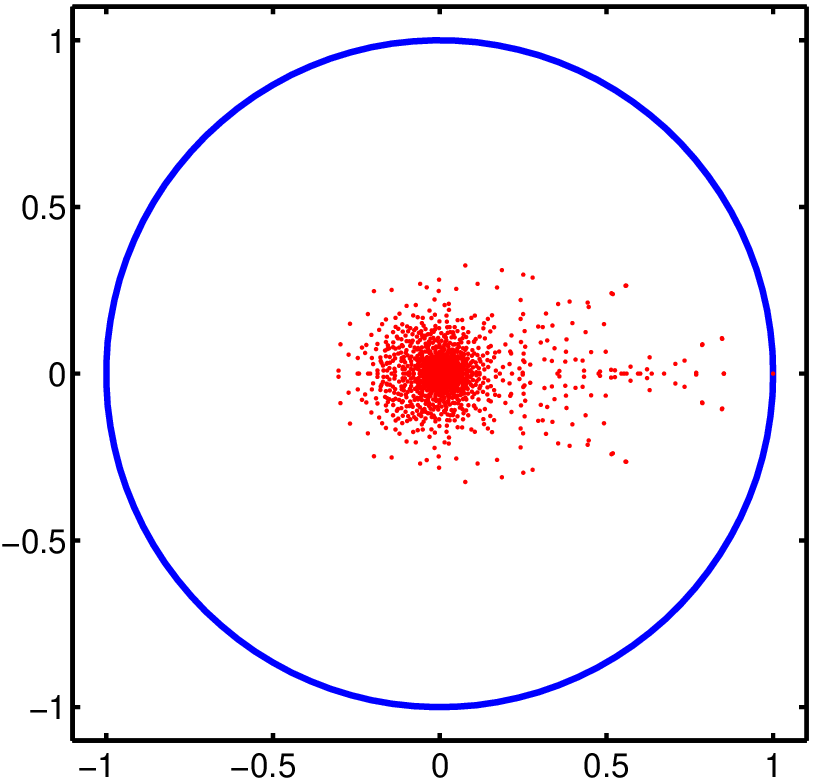}}
\caption{Plots of the eigenvalues of the multigrid iteration matrix
  for four different values of the relative stiffnesses.  The value of
  $\rho$ is the largest eigenvalue excluding the trivial eigenvalue of
  $1$.}
\label{spectrum:fig}
\end{figure}

In Figure \ref{spectrum:fig}, we plot the eigenvalues of the multigrid
iteration matrix in the complex plane for four different relative
stiffness.  We also report the spectral radius $\rho$ of the operator
on the space in which the mean pressure is set to zero.  Because the
pressure is unique only up to additive constants, there is always an
eigenvalue of $1$ that corresponds to pressure fields that are
constant on $\Omega_h$.  This trivial eigenspace does not affect the
convergence of the method.

The eigenvalues of the iteration matrix for the Stokes problem (i.e.,
with $\alpha = 0$) are shown in Figure \ref{spectrum:fig}(a).  All of
the nontrivial eigenvalues are clustered within a disc of radius
$0.24$ around the origin.  The spectrum is very similar at relative
stiffness 1.  In fact, the spectral radius is slightly smaller than in
the Stokes case.  At relative stiffness 10, most of the eigenvalues
are again clustered around the origin, but we now see a small number
of that are away from the origin.  The spectral radius is $0.30$,
which indicates a slight increase in the number of iterations over the
Stokes system ($\alpha = 0$).  We remark that we did not observe an
increase in the iteration count in our computational experiment until
the relative stiffness increased above 10.

At relative stiffness 100, the spectral radius is $0.85$, which is
consistent with our observation of slow convergence.  There are now a
notable number of ``large'' eigenvalues, but the majority of them are
still clustered near the origin.  For example, only about 5\% of the
eigenvalues are outside the disc of radius $0.25$, and only 1.5\% are
outside the disc of radius 0.5.  This eigenvalue distribution explains
why multigrid is a poor solver but an effective preconditioner.  It is
effectively damping a large eigenspace but poorly damping a small
eigenspace.  When used as a preconditioner for stiff problems,
multigrid acts to cluster most of the eigenvalues around $1$, and it
leaves a small set of scattered eigenspaces which are approximated by
a small Krylov space.  For a detailed analysis of this situation, see
Ref.\ \cite{Oosterlee:1998:EPM:277965.277972}.  Additionally, we note
the similarity to multiphase flow applications with sharp variation in
material properties.  In these applications, it has been observed that
multigrid is a poor solver but a very effective preconditioner for
Krylov methods \cite{Wright:2008:ERM:1405171.1405187,Sussman199981}.

\begin{figure}[hbt]
  \centering
  \subfigure[horizontal velocity]{\includegraphics[width=0.45\textwidth]{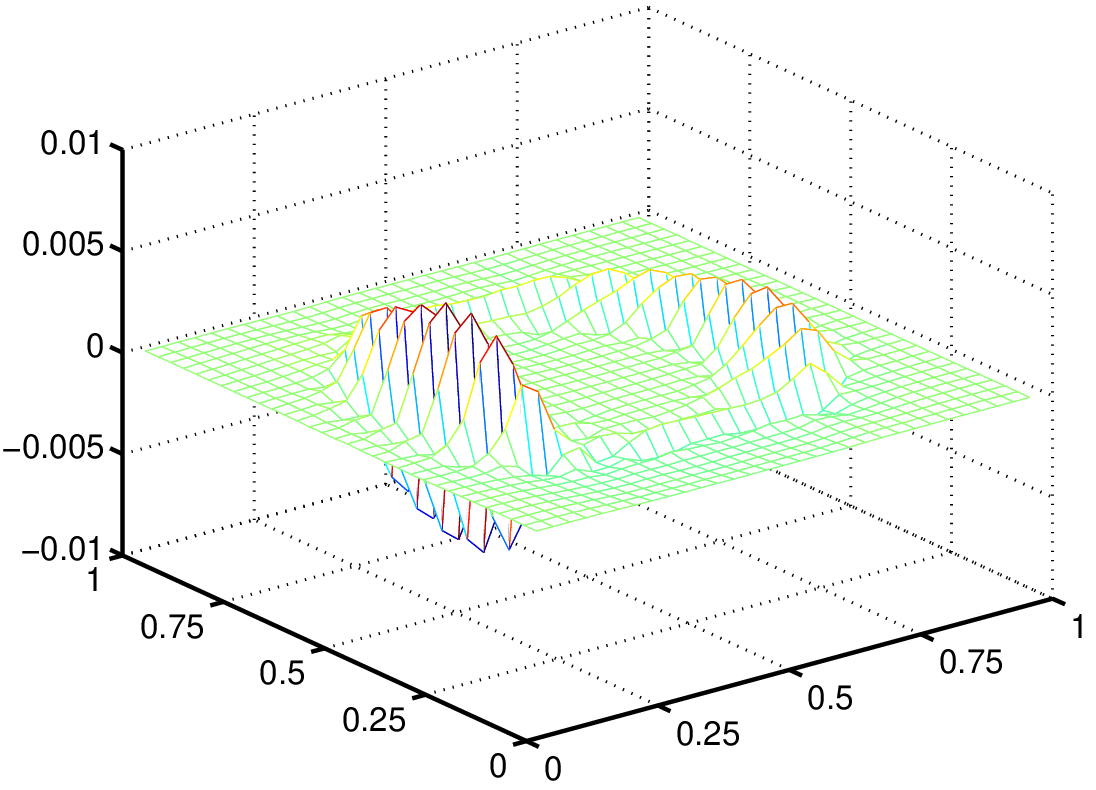}} \\
  \subfigure[vertical velocity]{\includegraphics[width=0.45\textwidth]{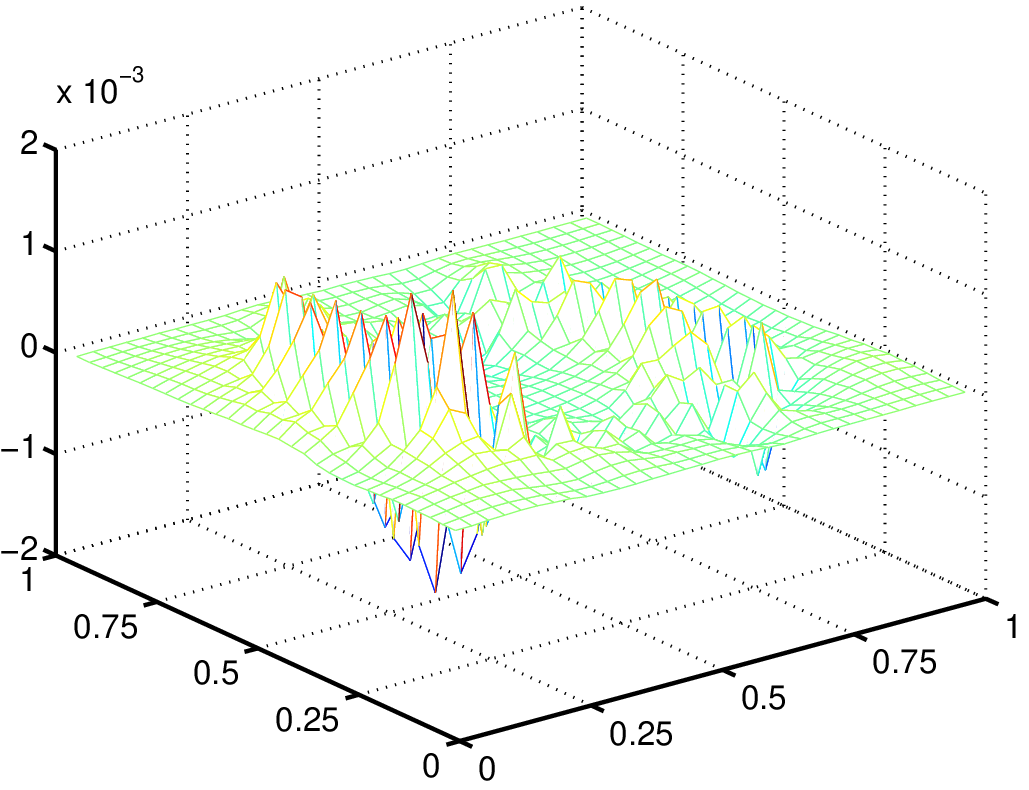}} \\
  \subfigure[pressure]{\includegraphics[width=0.45\textwidth]{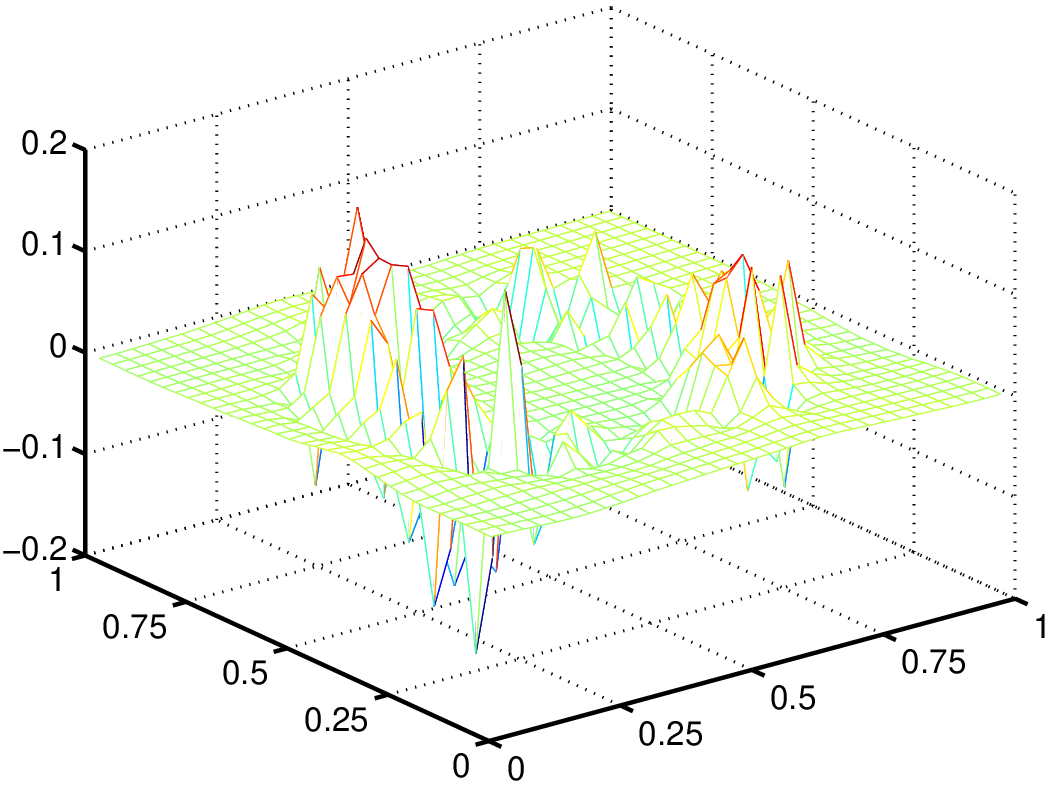}}
  \caption{Plots of the velocity and pressure corresponding to an
    eigenvector or the multigrid iteration matrix with eigenvalue
    0.8522 for relative stiffness $\alpha_\text{exp}/\alpha$.}
  \label{slow_mode:fig}
\end{figure}

We examine the velocity and pressure that correspond to the slowly
converging modes at large stiffness to attempt to obtain insight into
the poor performance of multigrid.  In Figure \ref{slow_mode:fig}, we
plot the two components of the velocity and pressure as functions of
space for the mode with eigenvalue $\approx 0.85$.  We see that this
mode exhibits high-frequency spatial oscillations, and that the
oscillations are concentrated near the immersed structure.  (Plots of
the other modes with large eigenvalues have a very similar features
and are not shown.)  This suggests that the poor performance of the
multigrid method results from our failure to smooth the high-frequency modes
associated with the elastic structure.

%
%
\section{Big-Box Smoothing}
\label{big-box_smoothing:sec}

In the previous section, we saw that for large elastic stiffnesses,
the multigrid solver with box relaxation required a large number of
iterations to converge or failed to converge because it did not smooth
high-frequency errors associated with the elasticity of the immersed
structure.  To understand this phenomenon better, we examine the
Eulerian elasticity operator $S_h\Af S_h^{*}$.

Recall that the operator $\Af$ represents the second-derivative
operator in the fiber direction (i.e., in the circumferential
direction in our annulus example), but with no cross-fiber coupling
(i.e., in the radial direction).  The interpolation operator $S_h^{*}$
maps Eulerian data to the Lagrangian mesh, and the spreading operator
$S_h$ maps Lagrangian data to the Eulerian grid.  Therefore the
operator $S_h \Af S_h^{*}$ is in some sense the projection of the
fiber-aligned Lagrangian second-derivative operator onto the Eulerian
grid.  The operator $\lf(\Delta_{h}+\alpha S\Af S^{*} \rt)$ thereby
resembles an anisotropic Laplacian operator, in which the degree of
the anisotropy increases with the elastic stiffness of the immersed
structure.

It is well known that strongly anisotropic problems require smoothers
that account for the anisotropy in order to obtain good multigrid
performance \cite{TROTTENBERG:BOOK}.  When the direction of anisotropy is
aligned with the grid, the standard approach is box-line smoothing, in
which all the cells in the $x$- or $y$-directions are collectively
updated.  For non-grid aligned anisotropy, such methods typically
alternate $x$-line and $y$-line smoothing in two spatial dimensions
\cite{TROTTENBERG:BOOK}.  The generalization of line smoothing (i.e.,
plane smoothing) to three spatial dimensions is computationally
expensive because it requires repeatedly solving two dimensional
problems.  Further, these approaches are more expensive than needed for
implicit IB methods because the anisotropy is localized to the region
covered by the immersed structure.

To improve the performance of our multigrid solver and preconditioner,
we follow a slightly different block relaxation approach that is still
in the spirit of the method introduced by Vanka \cite{Vanka86}.
Rather than update only the unknowns in a single cell, we instead
simultaneously update all unknowns in a rectangular box of size $n_x
\times n_y$, i.e., with $n_x$ cells in the $x$-direction and $n_y$
cells in the $y$-direction.  In the examples in this work, we take
$n_x=n_y$.  We call this approach \emph{big-box smoothing}.  As
before, the boxes are ordered lexicographically, and one smoothing
step involves sweeping over the boxes to update all of the velocities
and pressures associated with each $n_x \times n_y$ box.

Notice that the velocities that lie on the edges of the boxes are
updated twice, similar to the original Vanka scheme.  It is possible
to consider a further generalization of this scheme and consider
developing relaxation schemes with additional overlap between boxes.
We performed numerical experiments with different overlap sizes, and
we did not find a significant advantage in performance to justify the
added cost and complexity of including such additional overlap (data
not shown).

Our intent in developing such big-box smoothers is to provide a
relatively simple approach to smoothing oscillatory components arising
from the anisotropic coupling associated with the elasticity of the
immersed structure.  We note that the $x$- and $y$-line smoothers can
be considered as extreme cases of this scheme with $n_x=1$ and $n_y=N$
(or vice versa).  Unlike line smoothers, however, the size of the
boxes that we use does not change as the grid is refined.  Moreover,
our scheme naturally generalizes to three spatial dimensions.

\subsection{Solver performance}

We use the same test problem presented in the previous section.  We
begin with initial guess of zero for the velocity and pressure, and we
record the number of cycles needed to reduce the residual by a factor
of $10^{-6}$ for a range of stiffness.  The finest grid is
$32\times32$, and we use V-cycle multigrid with one presmoothing sweep
and one postsmoothing sweep ($\nu_1=\nu_2=1$).  We perform the test
for box sizes $n_x = n_y = 1$, $2$, $4$, $8$, and $16$.  \change{The
  coarsest grid in the V-cycle on which the exact solve is performed
  is $N_{c}\times N_{c}$ where $N_{c}=\max(4,n_{x})$.  In later tests
  the size of the coarsest grid is chosen in the same manor.}

\begin{figure}[th]
  \centering
  \subfigure[MG]{\includegraphics[width=0.48\textwidth]{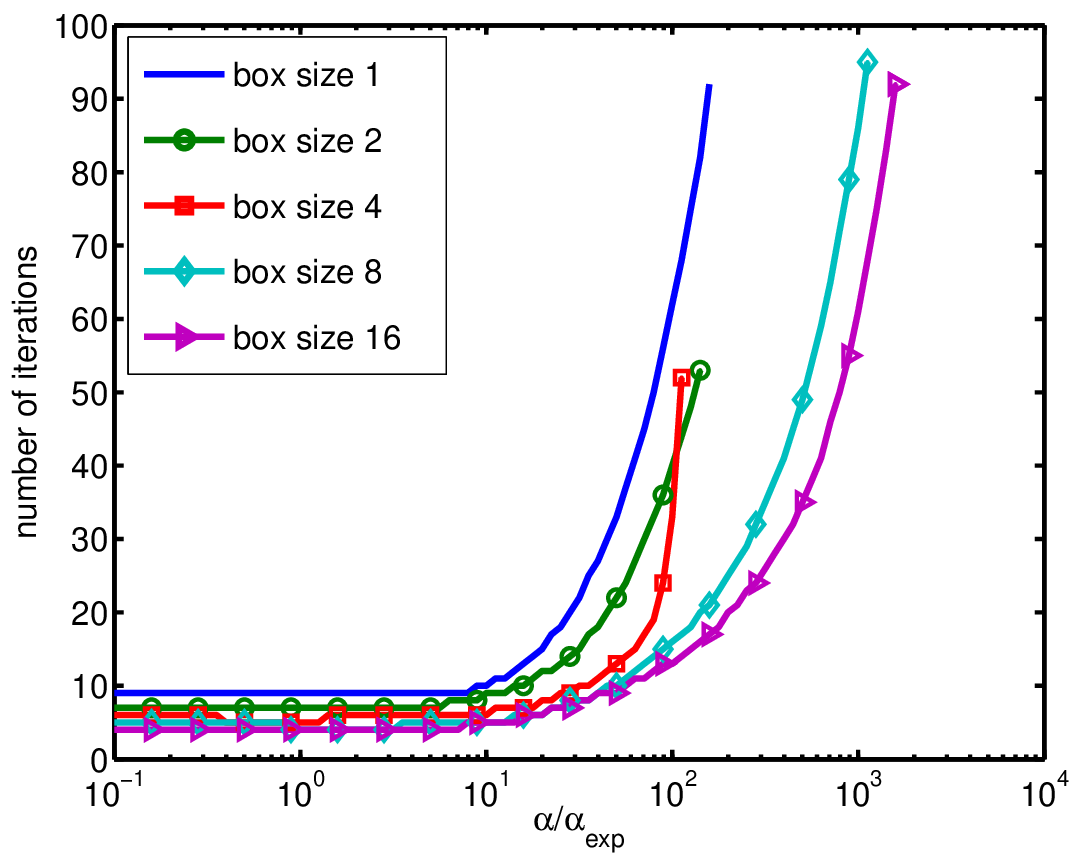}}
  \subfigure[MG-GMRES]{\includegraphics[width=0.48\textwidth]{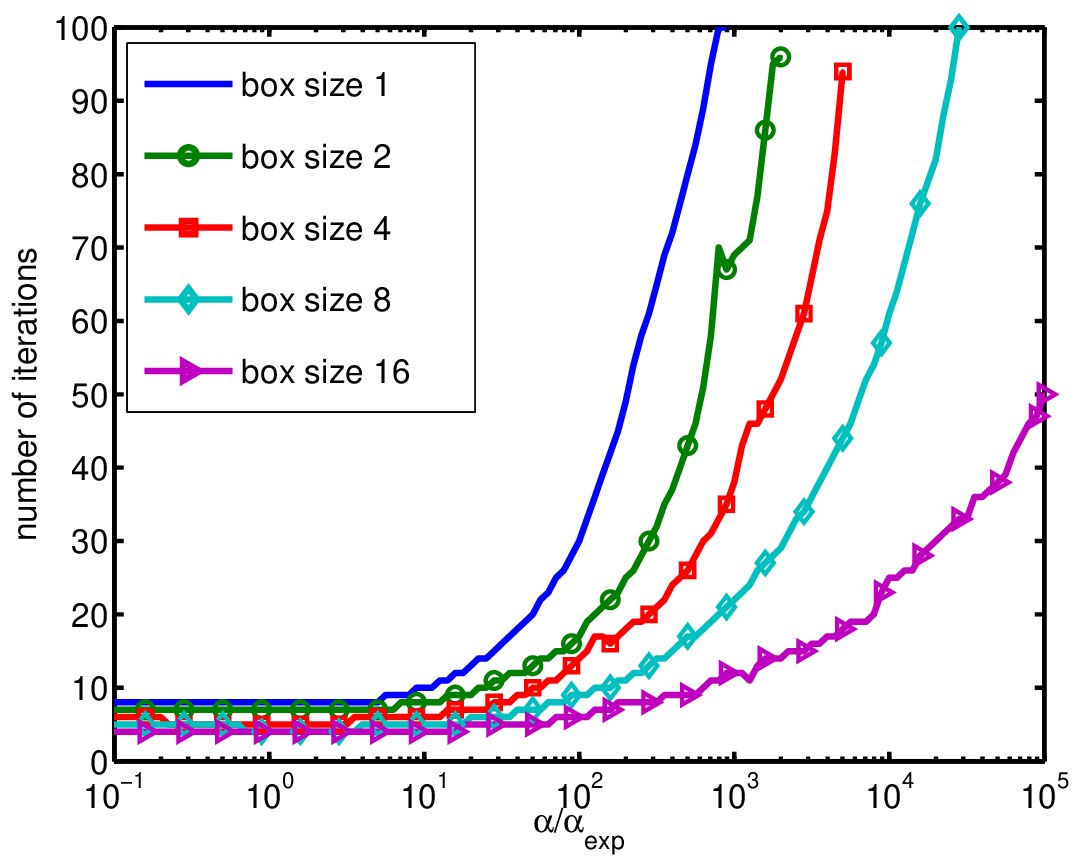}}
  \caption{Iteration counts for the multigrid and MG-GMRES methods to
    reduce the residual by a factor of $10^{-6}$ as a function of the
    relative stiffness for box smoothers of different sizes.}
   \label{big_box_cnt:fig}
\end{figure}

Figure \ref{big_box_cnt:fig}(a) shows iteration counts as a function
of stiffness for the stand-alone multigrid solver.  The single-box
smoothing results from the previous test are also included to
facilitate comparison.  Interestingly, the method with box sizes 2 and
4 fails at approximately the same stiffnesses as the single-box
smoother.  Before failure, the iteration counts are smaller with the
bigger boxes.  Box sizes 8 and 16 show significantly better
performance for stiff problems.  The iteration count does not start
increasing rapidly until a relative stiffness of about $100$.  The
method with box sizes 8 and 16 eventually fails, but at a stiffness
that is about an order of magnitude beyond the stiffnesses at which
the small-box smoothers fail.  Notice that the iteration counts with
box size 8 and box size 16 are very similar.

Figure \ref{big_box_cnt:fig}(b) shows the iteration counts for
MG-GMRES.  In general, the larger the box, the lower the iteration
count.  The difference in iteration counts between the box sizes is
particularly pronounced for relatively stiff problems.  As the box
size is increased, the iteration count rises less steeply as the
stiffness increases.  For example, for box sizes 2, 4, 8, and 16, the
iteration count increases from stiffness 0 to relative stiffness 1000
by a factors of 10, 6, 4, and 3, respectively.  These results indicate
that for very stiff problems, big-box smoothing in the multigrid
algorithm is a very effective preconditioner.  Because our tests
indicate that MG-GMRES is a much more robust method, from this point
on, we present test results only for MG-GMRES.


\subsection{Effect of grid refinement}

\begin{figure}[t]
 \centering
   \subfigure[$\alpha=6$; non-stiff]
           {\includegraphics[width=0.45\textwidth]{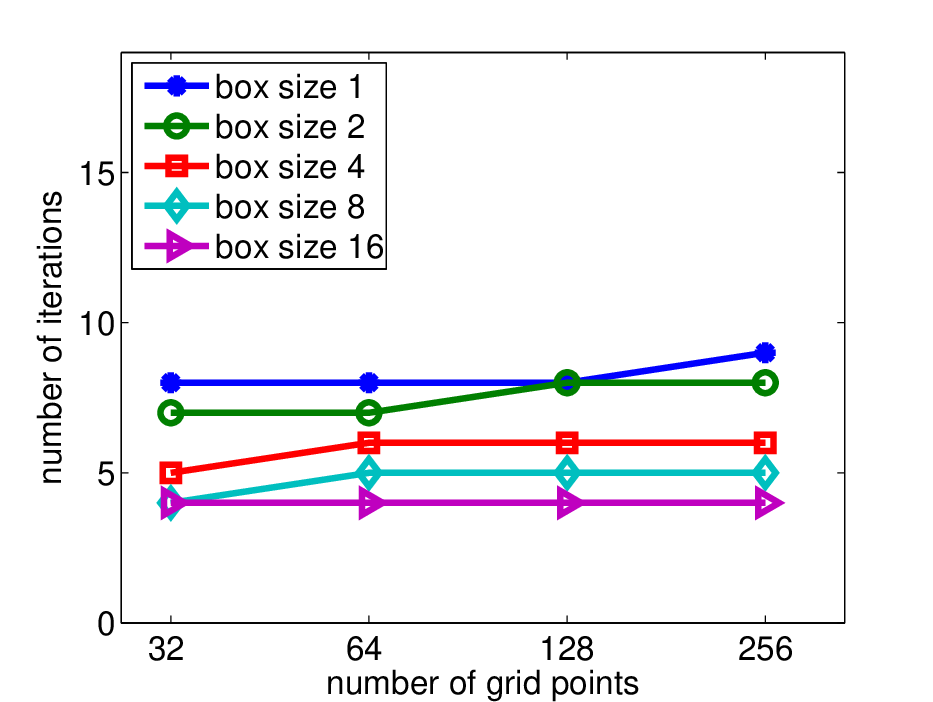}} \hfill
   \subfigure[$\alpha=60$; mildly stiff]
           {\includegraphics[width=0.45\textwidth]{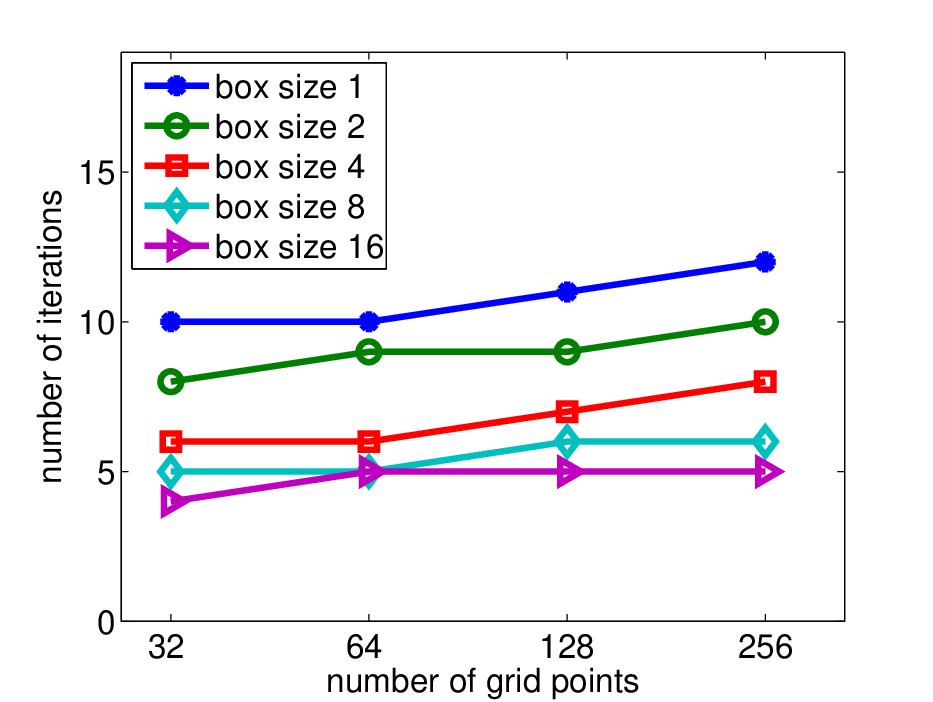}}
   \subfigure[$\alpha=600$; moderately stiff]
           {\includegraphics[width=0.45\textwidth]{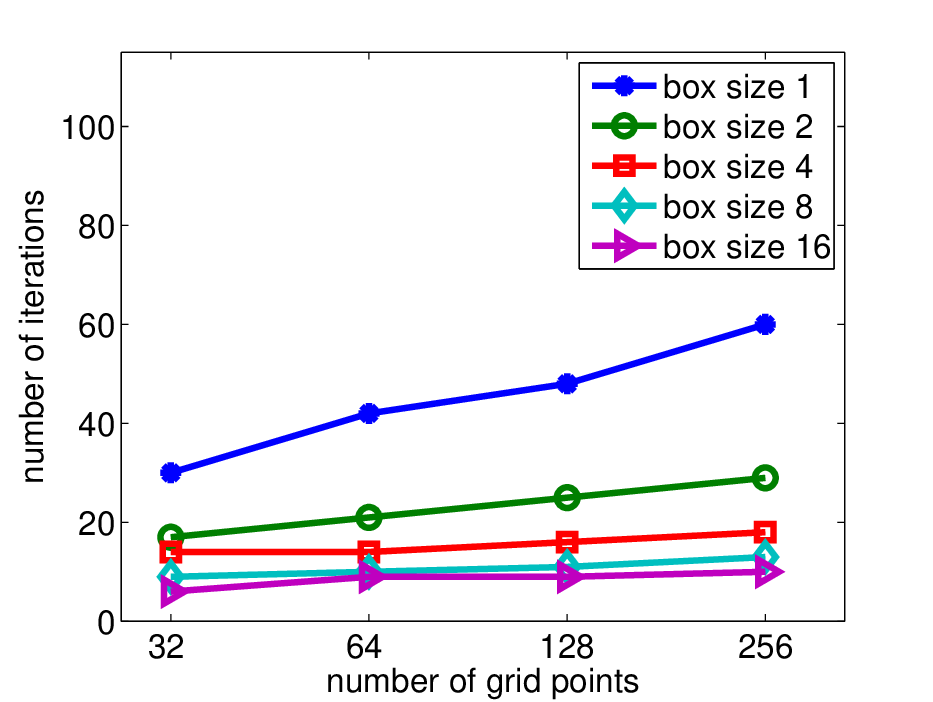}} \hfill
   \subfigure[$\alpha=3000$; very stiff]
           {\includegraphics[width=0.45\textwidth]{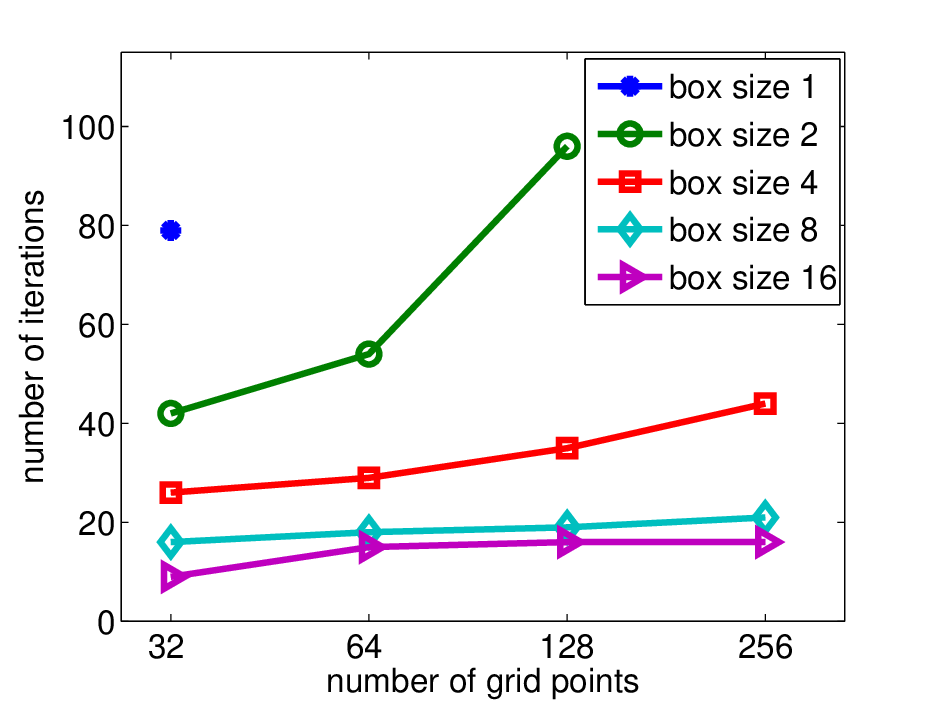}}
\caption{Iteration counts of MG-GMRES for different stiffnesses, box
  sizes, and different grid resolutions.  The number of grid points
  refers to the number of grid cells in each direction on the Eulerian
  grid.}
\label{gm_iter_cnt_grids:fig}
\end{figure}

So far, all of the tests were performed on an Eulerian grid with fixed
grid spacing (specifically, $N = 32$).  Here we explore the how grid
refinement affects the performance of the MG-GMRES algorithm in the
context of the same test problem used previously.  We choose four
values of the stiffness to explore: $\alpha=6$, $\alpha=60$,
$\alpha=600$, and $\alpha=3000$.  These values of the stiffness
roughly correspond to relative stiffnesses $1$, $10$, $100$, and $500$
on the $32\times32$ grid.  We chose these four values to characterize
regimes that are non-stiff, mildly stiff, moderately stiff, and very
stiff.

Figure \ref{gm_iter_cnt_grids:fig} shows the number of MG-GMRES
iterations needed to reduce the residual by a factor of $10^{-6}$ for
different box sizes and for four different grid resolutions.  As
before, the maximum number of iterations was set to 100.  For the
non-stiff and mildly stiff cases (Figure
\ref{gm_iter_cnt_grids:fig}(a,b)), the iteration count is essentially
independent of the grid size for all box sizes.  For the moderately
stiff and very stiff cases (Figure \ref{gm_iter_cnt_grids:fig}(c,d)),
the iteration count grows as the grid is refined for small boxes, but
the iteration count is essentially independent of the grid size for
the two largest box sizes (8 and 16).  For the very stiff case, the
smallest boxes often failed to converge in fewer than 100 iterations;
see data for box sizes 1 and 2 in Figure
\ref{gm_iter_cnt_grids:fig}(d).  These results are consistent with
those of the previous section.  For mildly stiff problems, small-box
smoothers perform well, but for stiff problems, small-box smoothers
require a large number of iterations to converge.  The difference in
performance between big-box smoothers and small-box smoothers is even
more striking as the grid is refined.


\subsection{Cost of smoothing and number of smoothing sweeps}

All of the previous results were generated with one presmoothing sweep
and one postsmoothing sweep ($\nu_1=\nu_2=1$).  Here we explore how
the number of smoothing sweeps affects the convergence of the MG-GMRES
algorithm for the total number of sweeps per level (i.e.,
$\nu_1+\nu_2$) ranging from 1 to 4.  If $\nu_1+\nu_2$ is even, then we
use an equal number of presmoothing and postsmoothing sweeps, and if
$\nu_1+\nu_2$ is odd, then we perform one additional presmoothing
sweep.  We consider box sizes $1$, $4$, and $8$ and relative
stiffnesses $10$, $100$, and $500$.  As before, these three different
stiffnesses characterize the mildly stiff, moderately stiff, and very
stiff regimes.  The test problem is the same as that used previously.
The Eulerian grid spacing is $h=2^{-6}$.

\begin{table}[h]
  \centering
  \caption{Iteration counts for different relative stiffnesses, box
    sizes ($b$), and numbers of presmoothing sweeps ($\nu_{1}$) and
    postsmoothing sweeps ($\nu_{2}$).  The total work, given in
    parentheses, is estimated as the number of iterations times
    $(\nu_1+\nu_2+1)$. For each box size and for each stiffness, the
    entry with the lowest total work is highlighted in bold.  (The
    number of MG-GMRES iterations was capped at 100, and work
    estimates are not provided for cases requiring more than 100
    iterations.)}
  \begin{tabular}{|c|c|c|||r|r|r|}
    \hline
    $b$ & $\nu_{1}$ & $\nu_{2}$ & $\alpha/\alpha_\text{exp}=10$ &    $\alpha/\alpha_\text{exp}=100$ &  $\alpha/\alpha_\text{exp}=500$ \\ \hline \hline
    1 & 1 & 0 &   16 (32)    &   53 (106)   &$\ge$100 (---) \\ \hline
    1 & 1 & 1 & {\bf 9 (27)} &   33 (99)    &      80 (240) \\ \hline
    1 & 2 & 1 &    7 (28)    & {\bf 21 (84)}& {\bf 55 (220)}\\ \hline
    1 & 2 & 2 &    6 (30)    &   17 (85)    &      44 (220 )\\ \hline \hline
    4 & 1 & 0 &   10 (20)    &   20 (40)    &   53 (106)    \\ \hline
    4 & 1 & 1 & {\bf 6 (18)} & {\bf 11 (33)}& {\bf 23 (69)} \\ \hline
    4 & 2 & 1 &    5 (20)    &    9 (36)    &   18 (72)     \\ \hline
    4 & 2 & 2 &    5 (25)    &    7 (35)    &   15 (75)     \\ \hline \hline
    8 & 1 & 0 &      8 (16)&       14 (28)&      28 (56)\\ \hline
    8 & 1 & 1 & {\bf 5 (15)}& {\bf 8 (24)}& {\bf 15 (45)}\\ \hline
    8 & 2 & 1 &      5 (20)&       7 (28)&       12 (48)\\ \hline
    8 & 2 & 2 &      4 (20)&       6 (30)&       10 (50)\\ \hline
  \end{tabular}
  \label{smooths:tab}
\end{table}

Table \ref{smooths:tab} shows the number of iterations needed to
reduce the residual by a factor of $10^{-6}$ along with a simple
estimate for the amount of computational work required to reach this
threshold.  The total estimated work to solve the problem is the
number of iterations times the work per iteration.  We estimate the
work per iteration as $(\nu_1+\nu_2+1)$; the ``plus one'' is included
to account for the work per iteration in addition to smoothing.  As
expected, the iteration count goes down as the number of smoothing
sweeps goes up.  For each box size and for each stiffness, the entry
in Table \ref{smooths:tab} with the lowest total work is highlighted
in bold.  For box sizes 4 and 8, one presmoothing sweep and one
postsmoothing is always the most efficient choice.  For the single-box
smoother, an additional smoothing sweep reduces the total work for
stiff problems.  As our previous results have demonstrated, for very
stiff problems, big-box smoothers are more effective.

We remark that our work estimate is based on the total number of
smoothing operations and is independent of the size of the box size.
This estimate is therefore useful only for comparisons in which the
box sized is kept fixed.  Because each step of the smoother requires
solving a linear system on each box, the cost of each smoother sweep
increases with box size.  However, as our results show, larger boxes
also reduce the total number of solver iterations required.  It seems
likely that there will not generally be a single set of algorithmic
options that is most efficient, but instead the optimal choices will
depend on the problem, implementation, and possibly the computer
architecture.

\section{Time-Dependent Problem}
\label{timedeptest:sec}

All of the previous tests focused on solving for the fluid velocity
for a prescribed, fixed structure position.  In this section, we test
the performance of MG-GMRES algorithm in a dynamic IB simulation.  We
place the same structure used in previous tests in a background shear
flow.  Unlike the previous tests, here the domain is a rectangle of
height 1 and length 2, and the background motion of the fluid is
driven by boundary conditions that, in the absence of the elastic
structure, would drive the shear flow $(u,v)=(y,0)$.  The structure is
initially centered at $\x_\text{c} = (0.5,0.5)$, and the simulation is
run until time $t=1$.

\begin{figure}[b]
  \centering
  \subfigure[$\gamma=60$]{\includegraphics[width=0.48\textwidth]{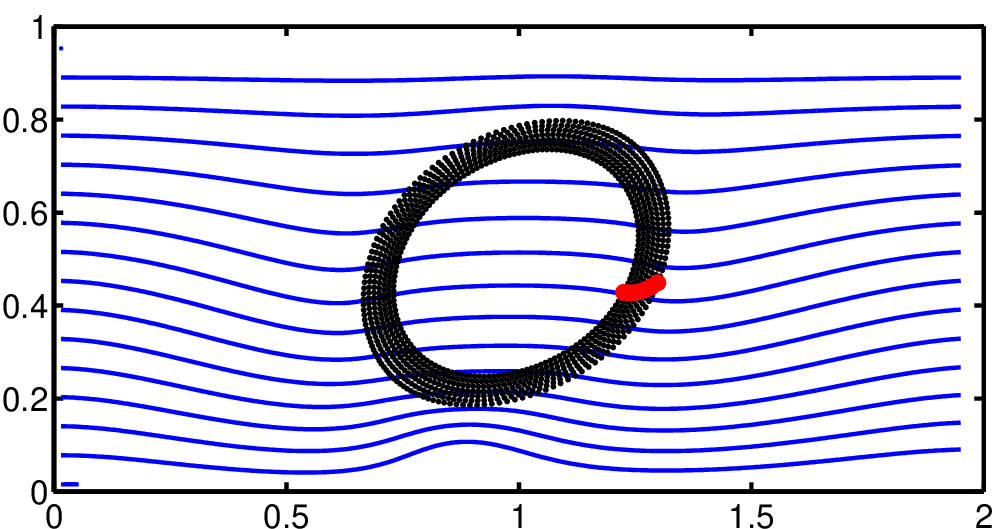}} \hfill
  \subfigure[$\gamma=240$]{\includegraphics[width=0.48\textwidth]{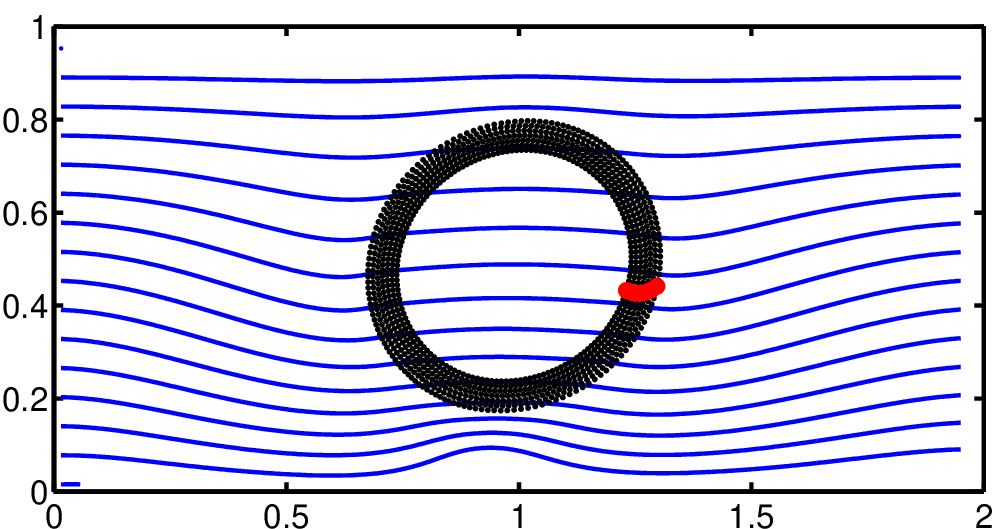}} \\
  \subfigure[$\gamma=2400$]{\includegraphics[width=0.48\textwidth]{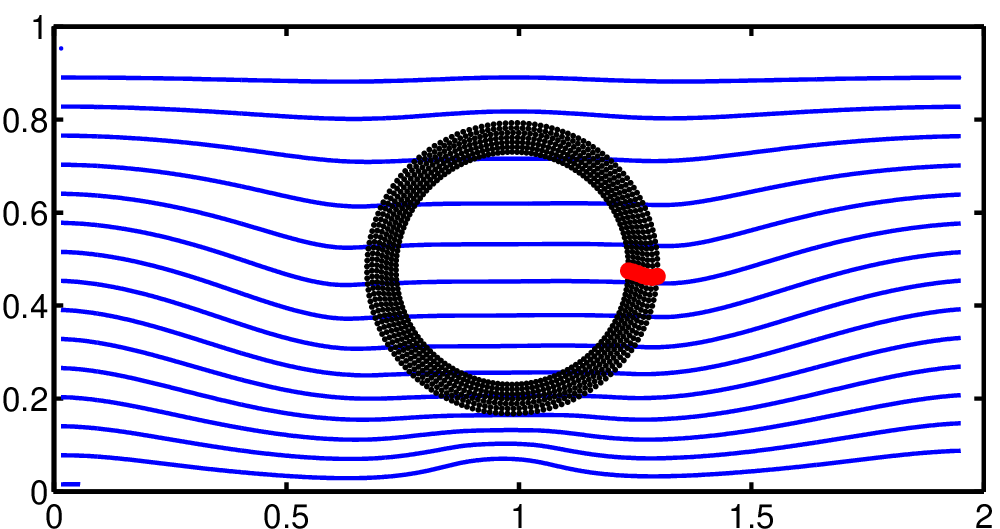}} \hfill
  \subfigure[$\gamma=24000$]{\includegraphics[width=0.48\textwidth]{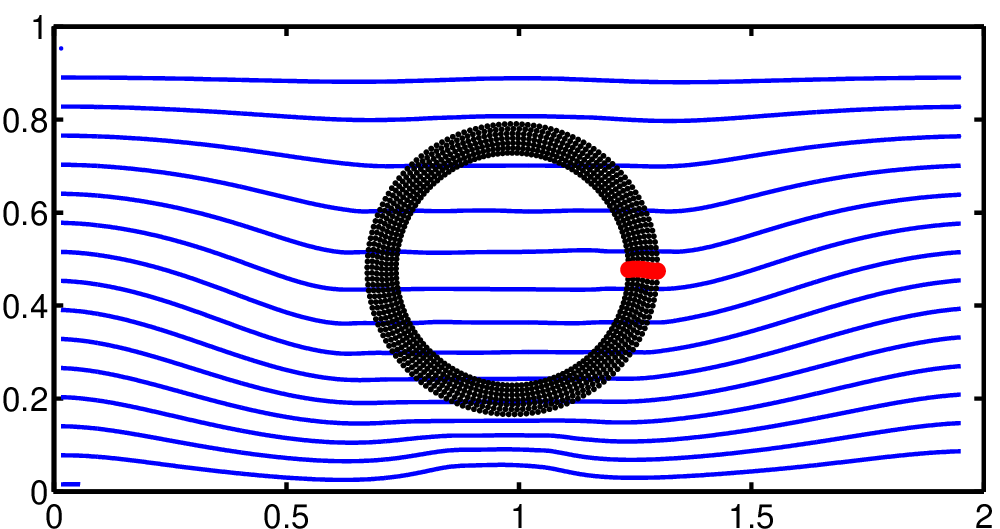}}
  \caption{Streamlines and structure location at time $t=1$ for a
    range of stiffnesses.  The red markers on the structure are used
    to highlight rotation and internal deformation.  At the beginning
    of the simulation the red markers were aligned horizontally in the
    $x$-direction.}
  \label{sheartest_flow:fig}
\end{figure}

We choose this test because the physical time scale is set by the
background flow, not by the stiffness of the structure.  This is the
type of problem for which implicit-time methods are needed.  In Figure
\ref{sheartest_flow:fig}, we show the structure's position at the end
of the simulation.  The elastic stiffness affects how much the
structure deforms, but the speed of translation is insensitive to the
stiffness.

\begin{figure}[ht]
  \centering
  \subfigure[]{\includegraphics[width=0.48\textwidth]{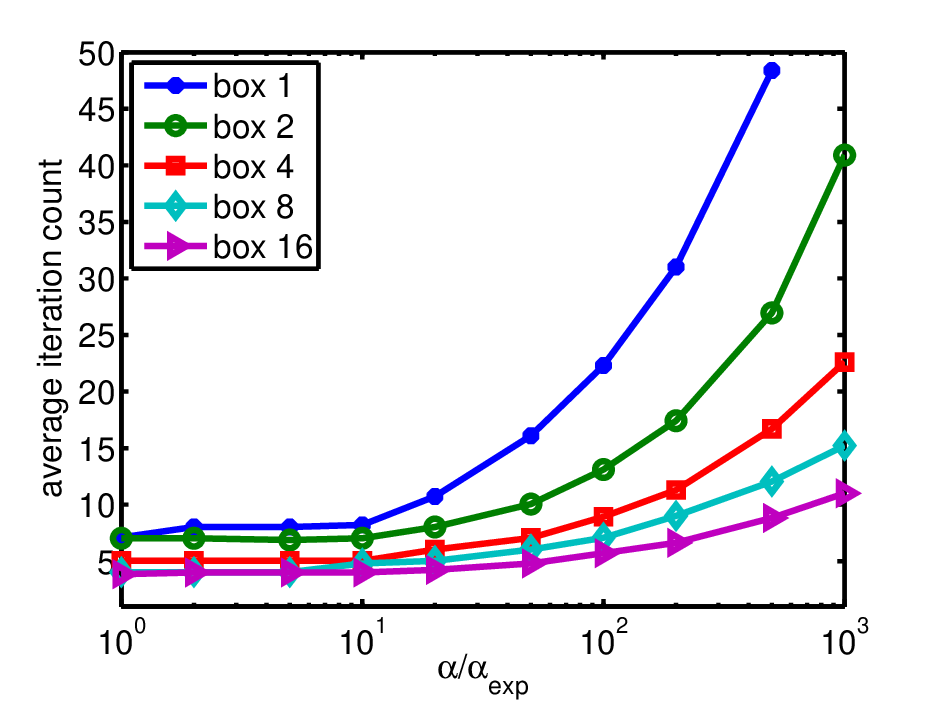}}
  \subfigure[]{\includegraphics[width=0.48\textwidth]{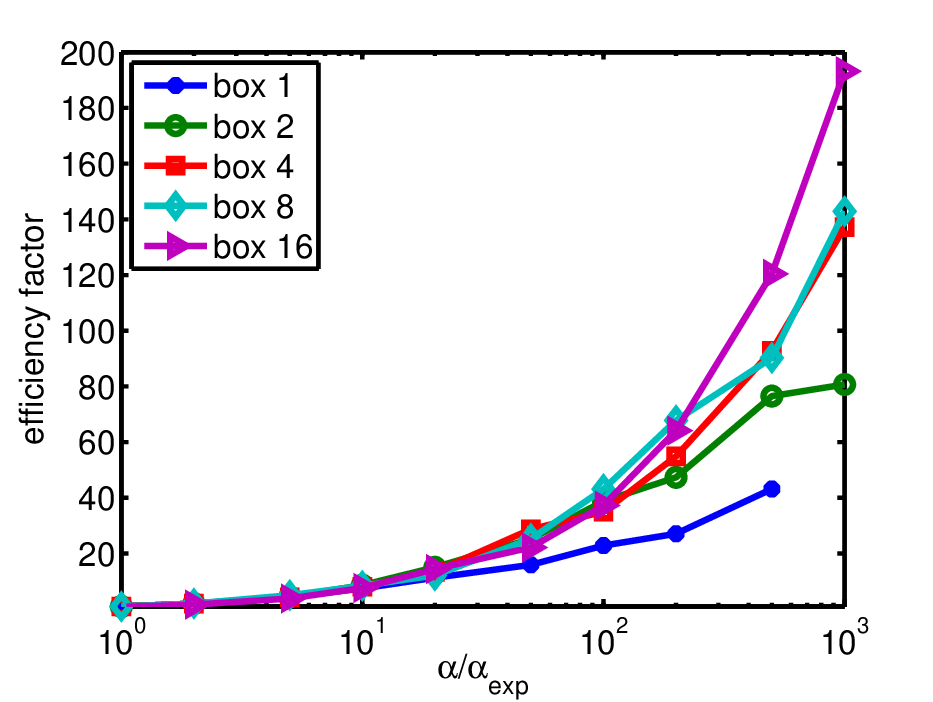}}
  \caption{(a) Average number of iterations of the MG-GMRES method per
    time step in the implicit-time simulation.  (b) The efficiency
    factor is the expected speed-up gained by using the implicit
    method in place of an explicit method.}
  \label{sheartest:fig}
\end{figure}

We discretize the Eulerian domain with grid spacing $h=2^{-5}$.  For
the implicit-time simulations, we fix the time step at $\dt=1/40$, and
we take $40$ time steps.  The maximum velocity is about $1$, and so
the Courant number of these simulations is about $0.8$.  Although
there is no stability constraint on the time step for this problem, we
choose to keep the Courant number less than one for reasons of
accuracy.  Recall that the time-dependent spreading an interpolating
operators are lagged in time, and so it is reasonable to require that
material points move less than a mesh width per time step.

We use the MG-GMRES method with $\nu_1=\nu_2=1$ to solve for the fluid
velocity to a relative tolerance of $10^{-6}$ at each time step.  Over
the course of a simulation, the number of iterations varies slightly
from time step to time step.  In Figure \ref{sheartest:fig}(a), we
report the average number of iterations of the solver per time step
over the course of the simulation as a function of the relative
stiffness for different box sizes.  As in the static tests, the
iteration counts are fairly constant up to a relative stiffness of
$10$.  Increasing the box size always lowers the iteration count, and
the difference in iteration count for different box sizes is
particularly striking for very stiff problems.

It is interesting to note that for all box sizes, the iteration count
grows sublinearly with the stiffness.  For example, at relative
stiffness 100, the iteration count increases by a factor of 1.5 (box
size 16) to 3 (box size 1) over the iteration count required to solve
the Stokes equations alone (i.e., without including the IB elasticity
operator).  This represents only a small increase in work compared to
an explicit-time method, which would require 100 times more time steps
to reach the same point in time.

We perform the same simulation using the explicit-time method to
compare the performance of the implicit and explicit method.  In the
explicit-time method, each time step involves solving the Stokes
equations for given IB forces.  We use the same MG-GMRES method with
big-box smoothing for solving the Stokes equations as we used for
solving the equations in the implicit-time method.  We use time step
sizes in the explicit-time simulations that are just below the
stability limit.  We perform 40 time steps and record the number of
MG-GMRES iterations.  We extrapolate to estimate the total number of
iterations to reach time $t=1$.

We define the efficiency factor as the ratio of the work of the
explicit-time method to the work of the implicit-time method needed to
complete this simulation.  The efficiency factor is the expected
speed-up one would gain by using the implicit method in place of the
explicit method.  To estimate the computational work, we use the
number of iterations of the MG-GMRES method.  In Figure
\ref{sheartest:fig}(b), we report the efficiency factor as a function
of the relative stiffness for different box sizes.  Up to a relative
stiffness of about 10, the efficiency factor is similar for all box
sizes.  As the stiffness increases, larger boxes outperform smaller
boxes.  These results show that for moderately stiff and very stiff
problems, the implicit method with big-box smoothing can be as much as
50--200 times more efficient than an explicit method.

\section{Discussion}
\label{discussion:sec}

The popularity of the immersed boundary method is driven by its
simplicity and robustness.  Implementations of the IB method that use
explicit-time solution algorithms generally require only solvers for
the fluid equations along with routines to compute elastic forces and
to transfer data between the Lagrangian mesh and the Eulerian grid.
The price of this simplicity is the severe restriction on the largest
stable time step permitted by such schemes.  One route to overcoming
this time step restriction is to develop implicit-time versions of the
IB method, but most previously developed implicit-time IB methods use
specialized algorithms to achieve substantial speed-ups over
explicit-time methods.

The goal of this work is to investigate solution approaches to
implicit IB methods that balance efficiency, robustness, and
simplicity.  A distinguishing feature of our method is that we
formulate the problem on the Cartesian grid.  This formulation allows
the use of standard tools from geometric multigrid to solve the
equations.  Our algorithm is similar to multigrid methods for the
Stokes equations that use coupled smoothing and is relatively
straightforward to implement.  The major differences between the
present algorithm and standard multigrid methods for incompressible
flow are the presence of IB elasticity operator $S_h \Af S_h^{*}$, and
the use of big-box smoothing.  To our knowledge, big-box smoothing has
not been proposed previously in the context of incompressible flow
problems.  Moreover, multigrid-based approaches to developing solvers
for implicit IB methods that do not eliminate the Lagrangian degrees
of freedom would generally need to resort to approaches such as
algebraic multigrid to generate coarse-grid versions of similar
operators for cases in which the immersed structure has complex
geometry, as is often the case in practice.  In our approach, because
the relevant elasticity operator is defined on the Cartesian grid, it
is straightforward to construct coarse-grid versions of this operator
using Galerkin coarsening using the structured-grid restriction and
prolongation operators.

Potential limitations of this study are that we restricted our
exploration to zero Reynolds number flow and to structures with
nonzero thickness with linear constitutive laws.  Extensions to
nonzero Reynolds numbers are straightforward.  In fact, solving the
time-independent problem, as we do herein, is \emph{more} challenging
than solving the time-dependent problem.  \change{The time-dependent
  problem presents other possible solution methods and preconditioning 
  strategies based on fractional stepping techniques, which do not
  extend to the time-independent problem.}  Provided that the
nonlinear convection terms are treated explicitly, the system of
equations that must be solved at each time step is very similar to
\eqref{matrixeq:eq}, with the main difference being that the operator
applied to the velocity in the time-dependent problem is
$I-\dt\Delta_{h}-\dt^2\gamma S_h\Af S_h^{*}$.  The presence of the
identity matrix and the appearance of an additional factor of $\dt$
both result in a better conditioned system, and likely will improve
the convergence of our method.

In many applications of interest, the elastic structure is modeled as
a thin interface of zero thickness.  Our method applies equally well
to this problem.  We choose to use ``thick'' structures for our
numerical tests because the IB method is known to yield poor volume
conservation when applied to problems involving thin elastic
structures \cite{GRIFFITH:2012:IBVOLCONSERVE}.  This unphysical
leakage of fluid across the boundary is exacerbated by high elastic
stiffness, making it difficult to reach the extreme stiffnesses
considered in this work.  (Although not shown here, we performed a
series of tests with thin boundaries, but the leakage limited our
tests to moderate stiffnesses.)  Our multigrid solution algorithm
performs equally well on thin and thick boundaries, however.  The
major difference between these two arose upon grid refinement.  For a
fixed elastic stiffness, thin boundaries become more numerically stiff
as the grid is refined.  That is, the time step restriction scales
like the grid spacing.  If the time step is reduced simultaneously
with the grid spacing, as is necessary in practice to maintain a fixed
Courant number, then the solver performance is essentially grid
independent.

\change{In the test problems considered in this paper, the structure
  inherits the viscosity as the underlying fluid.  Immersed boundary
  methods which include a separate structure viscosity have been
  developed
  \cite{FAI:2013:VarVisc,Huang20092650,STRYCHALSKI:2012:VEIB}. A
  semi-implicit-time discretization which accounts for structure
  viscosity was proposed and analyzed in \cite{STRYCHALSKI:2012:VEIB}.
  The resulting equations are very similar to those we present in this
  manuscript, the difference being the form the boundary force
  operator $\Af$, and so our algorithm extends naturally to this
  type of problem.  }

Finally, nonlinear constitutive laws could be treated by semi-implicit
time stepping schemes that effectively linearize the nonlinear
equations.  Alternatively, nonlinear time discretizations could be
solved via Newton's method.  In either case, it is necessary to solve
repeatedly linearized problems of the form considered herein.

A possible criticism of our method is that it requires the use of
``big'' boxes to achieve a robust solution method.  A potential
concern with larger boxes is the computational expense of solving the
restriction of the full equations to the boxes.  We remark that in $d$
spatial dimensions, each $n \times n$ box has $d (n+1) n^{d-1}$
velocity components and $n^d$ pressure components.  For $d=2$ and
$n=4$, this is a total of 56 degrees of freedom, and if we were to use
a dense representation of the box operator, the memory requirements
would be approximately 25 KB, which is small enough to fit into
high-speed L1 cache on most modern CPUs.  In this case, the cost of a
direct solver for the box operator will be essentially negligible.
For $d=3$ and $n=4$, each box has 304 degrees of freedom, and the
corresponding dense representation requires approximately 750 KB,
which is too large to fit into L1 cache, but which fits comfortably in
the L2 or L3 cache available on most systems.  Although our tests show
that a box size of four does not yield perfect grid-independent
convergence rates for extremely stiff systems, it does converge in a
reasonable number of iterations and yield good performance compared to
an explicit method.  We expect that for $n=4$, the local solves will
be efficient in a high-quality implementation of the algorithm in
either two or three spatial dimensions.  Of course, there may be cases
in which box sizes greater than four are needed.  In these cases,
implementations may also be able to exploit the sparsity of the box
operators.  (Only the ``(1,1) block'' of the box operator will be
relatively dense.)  Further, since the restriction of the full
operator to the boxes is itself a discrete saddle-point system, for
sufficiently large boxes, it may be worthwhile to investigate
approximate solution methods.

While our algorithm is simple in spirit, as the foregoing discussion
suggests, producing an optimized implementation is not trivial.  We
estimate that our algorithm will achieve a substantial speed-up over
explicit-time methods, but our efficiency estimates are based on
iteration counts for a fixed box size, and we did not compare the
efficiency across different box sizes. Our initial implementation of
the algorithm was not designed with wall-clock efficiency in mind, and
fully quantifying the performance of the algorithm requires that we
develop an optimized code.  

\change{An efficient explicit-time method could exploit that the
  subsystems associated with the big-box smoother are identical (away
  from the domain boundaries) and could be pre-factored and stored.
  Near the immersed boundary the matrix $S_h \Af S_h^{*}$ must be reformed at
  each time step, and the subsystems are all distinct.  Reforming $S_h
  \Af S_h^{*}$ does not present a high computational cost because each
  of these matrices is a sparse matrix with the number of elements
  proportional the number of immersed boundary points, but this is an
  additional cost the implicit method that is not reflected in
  iteration counts alone.  } However, given the differences in total
numbers of iterations required by the implicit and explicit schemes,
we believe that it is reasonable to expect to obtain at least
order-of-magnitude speed-ups over explicit-time solvers by using
optimized implementations of the present implicit-time algorithm.  We
are in the process of developing optimized versions of our algorithm
and including them in the IBAMR software package
\cite{IBAMR-web-page}, and we hope to report results from applications
of this code in the future.



\bibliographystyle{spmpsci}      
%

\end{document}